\documentclass[11pt]{article}      
\usepackage{latexsym,amsmath,amsfonts,amscd,amssymb}
\usepackage{amsthm}
\usepackage{amssymb,multirow}
\usepackage{graphics}
\usepackage{epsfig}
\usepackage{subfigure}
\usepackage{verbatim}
\usepackage{epstopdf}
\usepackage{color,leftidx}

\providecommand{\keywords}[1]{\textbf{Keywords:} #1}

\newtheorem{thm}{Theorem}[section]
\newtheorem{rem}{Remark}[section]

\newtheorem{Algorithm}{Algorithm}[section]
\newtheorem{example}{Example}[section]

\newcommand{\be}{\begin{equation}}
\newcommand{\ee}{\end{equation}}

\newcommand{\ben}{\begin{equation*}}
\newcommand{\een}{\end{equation*}}

\author{Zheng Chen\footnote{Department of Mathematics, Iowa State University, Ames, IA 50011. Email: zchen@iastate.edu}, Hongying Huang\footnote{School of Mathematics, Physics and Information Science, Zhejiang Ocean University; Zhejiang and Key Laboratory of Oceanographic Big Data Mining \& Application of Zhejiang Province, zhoushan, Zhejiang, China, Email: huanghy@lsec.cc.ac.cn},
Jue Yan\footnote{Department of Mathematics, Iowa State University, Ames, IA 50011. Email: jyan@iastate.edu}}
\title{Third order Maximum-Principle-Satisfying Direct discontinuous Galerkin methods for time dependent convection diffusion equations on unstructured triangle mesh}
\date{\today}

\begin{document}

\maketitle

\begin{abstract}
We develop 3rd order maximum-principle-satisfying direct discontinuous Galerkin methods
\cite{liu2009direct, liu2010direct, VY-DDG-2010, yan2013new}
for convection diffusion equations on unstructured triangular mesh. We carefully calculate the normal derivative numerical flux across element edges and prove that, with proper choice of parameter pair $(\beta_0,\beta_1)$ in the numerical flux, the quadratic polynomial solution satisfies strict maximum principle. The polynomial solution is bounded within the given range and third order accuracy is maintained. There is no geometric restriction on the meshes and obtuse triangles are allowed in the partition. A sequence of numerical examples are carried out to demonstrate the accuracy and capability of the maximum-principle-satisfying limiter.
\end{abstract}


\section{Introduction}
\numberwithin{equation}{section}

In this article, we study direct discontinuous Galerkin  finite element method \cite{liu2009direct} and its variations \cite{liu2010direct, VY-DDG-2010, yan2013new} to solve two-dimensional convection diffusion equations of the form,
\be{\label{eqn:2DConvecDiff}}
u_t  + \nabla \cdot F(u) - \nabla \cdot (A(u)\nabla u) = 0, \qquad  (x,y,t) \in \Omega\times(0,T),
\ee
with zero or periodic boundary conditions. We have spacial domain $\Omega \subset \mathbb{R}^2$ and initial condition $u(x,y,0) = u_0(x,y)$. The convection flux is denoted as
$F(u) = (f(u), g(u))$ and diffusion matrix $ A(u)=(a_{ij}(u))$ is assumed symmetric and positive definite.

On the continuous level, solution of (\ref{eqn:2DConvecDiff}) may satisfy the maximum principle, which states the evolution solution $u(x,y,t)$ being bounded below and above by the given constants,
$m\leq u(x,y,t) \leq M$. Here $m$ and $M$ are the lower and upper bounds of the initial and boundary data. It is desirable that the numerical solution satisfies the \emph{discrete maximum principle}.
The discrete maximum principle can be considered as a strong $L^{\infty}$ sense stability result. Failure of preserving the bounds or maintaining the positivity of the numerical solution may lead to ill-posed problems and practically cause the computations to blow up.
Thus it is attractive to have the numerical solution satisfy discrete maximum principle (or preserve positivity). Solution of equation (\ref{eqn:2DConvecDiff}) may represent a specific physical meaning and is supposed to be positive, thus negative value approximation loses physical meanings in such cases.

Generally it is very difficult to design high order numerical methods that satisfy discrete maximum principle for convection diffusion equations (\ref{eqn:2DConvecDiff}). No finite difference method is known to achieve better than second-order accuracy \cite{Evje2000, zhang-MPS-Review, XuZhengFu2014} that satisfies discrete maximum principle.  Much less is known for higher order methods such as spectral FEM, hp-FEM or finite volume methods \cite{Bertolazzi-Manzini2005, Guo2008, Sheng-Yuan2011, MPS-FV}. Compared to the elliptic type, more restrictive conditions on mesh are required to obtain discrete maximum principle for the parabolic type equations, see \cite{Vejchods2004, ThomeeBook2006, Holst2008, Thomee-MP2008}.

In this article, we study direct discontinuous Galerkin  method \cite{liu2009direct} and its variations {\cite{liu2010direct, VY-DDG-2010, yan2013new}}, and prove the polynomial solution satisfy discrete maximum principle with third order of accuracy. Discontinuous Galerkin (DG) method is a class of finite element method that use completely discontinuous piecewise functions as numerical approximations. Since the basis functions can be completely discontinuous, these methods have the flexibility that is not shared by standard finite element methods, such as the allowance of arbitrary triangulations with hanging nodes, complete freedom of choosing polynomial degrees in each element ($p$ adaptivity), and extremely local data structure and the resulting  high parallel efficiency.

Recently in \cite{zhang2010maximum, zhang-MPS-Review, zhangxxTri}, Zhang and Shu designed a maximum-principle-satisfying limiter for high order DG and finite volume methods for hyperbolic conservation laws. The key step in Zhang and Shu's discussion is to show the polynomial solution average falling in the given minimum and maximum bounds. For hyperbolic type equations, the solution average evolution only relies on the solution polynomial values on the element edges. For diffusion type equations, the evolution of solution average depends on the {\emph{solution derivative}} values on the edges, thus the technique developed in \cite{zhang2010maximum} can not be applied.

In \cite{liu2009direct}, we developed the direct DG method (DDG) as a new diffusion solver. The key contribution of direct DG method is the introduction of numerical flux to approximate the solution derivative at the discontinuous element boundaries. The scheme is directly based on the weak formulation of diffusion equation, thus gains its name the direct DG method. Now let's use the simple 2-D heat equation to go through the main idea of direct DG method,
\be\label{eqn:2dheat}
u_t-\bigtriangleup u=0.
\ee
Multiply the heat equation with test function $v$, integrate over the element $K$, have integration by parts and formally we obtain,
\begin{equation*}
\int_K u_t v~dxdy  - \int_{\partial K} \widehat{u_{\textbf{n}}} v~ds + \int_K \nabla u \cdot \nabla v ~dxdy= 0.
\end{equation*}
The numerical flux $\widehat{u_{\textbf{n}}}$ introduced in \cite{liu2009direct} is defined as follows,
\begin{equation*}
\widehat {u_{\textbf{n}}} =\widehat {\nabla u\cdot \textbf{n}}=\beta_{0}\frac{[u]}{h_K}+ \overline{\frac{\partial u}{\partial \textbf{n}}}
+\beta_{1} h_K [u_{\textbf{nn}}].
\end{equation*}
It involves the jump, the derivative average and higher order derivative jumps to approximate the normal derivative $u_{\textbf{n}}$ on the element boundary $\partial K$.
Here $\textbf{n}=(n_{1},n_{2})$ is the outward unit normal along ${\partial K}$ and
$h_K$ is the diameter of element $K$. The coefficient pair $(\beta_0, \beta_1)$ is chosen to guarantee the convergence of the scheme.

Due to accuracy loss of the original DDG method \cite{liu2009direct}, we further developed DDG method with interface correction in \cite{liu2010direct} in which optimal $(k+1)$th order convergence is obtained with any order $P^k$ polynomial approximations. We also have the symmetric \cite{VY-DDG-2010} and nonsymmetric version \cite{yan2013new} of the DDG methods. In this paper, we mainly carry out the maximum principle study on DDG method with interface correction \cite{liu2010direct} since it is the most efficient solver for time dependent diffusion equations. The maximum principle arguments discussed in the following sections also apply to DDG method \cite{liu2009direct} and its symmetric and nonsymmetric variations \cite{VY-DDG-2010, yan2013new}.

In \cite{Yan-MPS}, we prove the DDG solutions satisfy discrete maximum principle on rectangular and uniform triangular meshes with 3rd order of accuracy. We use an algebraic methodology and a monotonicity argument to show the polynomial solution average being bounded within the given range. The DG  polynomial solution was written out in the Lagrange format with the unknowns carefully chosen on the element. With Euler forward in time, we show the solution average at next time level depends on the current time level solution values in a monotone fashion. For unstructured mesh with possible obtuse triangles, it is very hard to identify six degrees of freedom to represent the $P^2$ quadratic polynomial solution such that the monotone argument in \cite{Yan-MPS} can be applied.

In this article we extend maximum principle studies of (\ref{eqn:2DConvecDiff}) on unstructured triangle mesh. Again let's use the heat equation to illustrate the new technique to carry out the proof. Notice that the key step of the discussion is to show the solution average falling in the given range. Take test function $v=1$ in the DDG scheme and discretize in time with Euler forward, we have the solution average evolving in time as,
\begin{equation*}
\label{2DAverageEvo-1}
\overline{u}^{n+1}_K =\overline{u}^{n}_K+\frac{\Delta t}{area(K)}\int_{\partial K} \widehat{u_{\textbf{n}}}~ds,
\end{equation*}
with the average defined as $\overline{u}^{n}_K=\frac{1}{area(K)}\int_K u^n_K(x,y)~dxdy$ and $u_{K}^{n}(x,y)$ as the polynomial solution at time step $t_n$ in element $K$.

Instead of identifying suitable locations as degrees of freedom and writing out $u_{K}^{n}(x,y)$ in the Lagrange format as in \cite{Yan-MPS}, we directly calculate the normal derivative flux $\widehat{u_{\textbf{n}}}$ from
the given solution values in element $K$ and its neighbors. Given suitable choice of coefficient pair $(\beta_0,\beta_1)$ in the numerical flux, we can bound the solution average  $\overline{u}_{K}^{n+1}\in [m,M]$ at time level $t_{n+1}$ once we know $u_{K}^{n}(x,y)\in [m,M]$ at previous time level $t_n$. Finally we borrow the maximum principle discussion of \cite{zhang-MPS-Review} to show the DG polynomial solution of general convection diffusion equations (\ref{eqn:2DConvecDiff}) satisfy strict maximum principle with 3rd order of accuracy. A sequence of numerical examples are carried out to demonstrate the DG solutions are strictly bounded by the given values and at the same time maintain the 3rd order accuracy. Solutions to nonlinear porous medium equations with nonnegative initial data are maintained sharply nonnegative. Examples of incompressible Navier-Stokes equations with high Reynold numbers are tested. Overshoot and undershoots are removed with maximum principle limiter applied.

The key feature of direct DG methods is the introduction of numerical flux $\widehat{u_{\textbf{n}}}$ that approximates the solution derivative $u_{\textbf{n}}$ on the discontinuous element boundary. This gives direct DG methods the extra flexibility and advantage over IPDG method \cite{Arnold1982} and LDG method \cite{Cockburn-Shu-1998}. Following this maximum principle framework, both IPDG method and LDG method can be proved to satisfy maximum principle with up to second order of accuracy. In \cite{zhang-Shu-convdiff}, DG solutions  with piecewise linear polynomial approximations are shown satisfying maximum principle on unstructured triangle mesh.

The paper is organized as follows. We first review the scheme formulation of direct DG method with interface correction \cite{liu2010direct} in section \ref{sec:DDGs}. In Section \ref{sec:3rdMPSDDG}, we prove the direct DG solutions satisfy discrete maximum principle with 3rd order accuracy. We conduct numerical tests to validate the theoretical results in Section \ref{sec:numerical}. Section \ref{sec:appendix} serves as the Appendix in which we provide one way to construct a specific quadrature rule for quadratic polynomials  with selected points as quadrature points.

\numberwithin{equation}{section}
\section{Direct DG method with interface correction}\label{sec:DDGs}

We first recall the scheme formulation of direct DG method with interface correction \cite{liu2010direct} for two-dimensional diffusion equations,
\be{\label{eqn:2DDiff}}
u_t  - \nabla \cdot ( A(u)\nabla u) = 0, \qquad  (x,y,t) \in \Omega\times(0,T),
\ee
with initial condition $u(x,y,0) = u_0(x,y)$ and zero or periodic boundary conditions. The complete scheme formulation of convection diffusion equation (\ref{eqn:2DConvecDiff}) will be laid out toward the end of this section. We should specify the DG method is for spatial discretization and we will incorporate high order TVD Runge-Kutta methods \cite{Shu-Osher-1989,Shu-Osher-TVD-1988} to march forward the solution in time. As an explicit scheme, our method is thus more efficient for convection dominated problems. However, the extremely local dependency allows a very efficient parallelization and dramatically improves the efficiency of the explicit method.

Let $\mathcal{T}_h$ be a shape-regular partition of the polygonal domain $\Omega$ into triangle elements $\{K\}_{K\in \mathcal{T}_h}$ with $\overline{\Omega}=\cup_{{K\in \mathcal{T}_h}}\overline K$. By $h_K=\text{diam}(K)$, we denote the diameter of the triangle element $K\in \mathcal{T}_h$. We denote  $h=\max_{K\in \mathcal{T}_h}h_K$ as the mesh size of the partition.
We have $P^k(K)$ representing the  $k$th degree polynomial space on element $K$. The DG solution space is defined as,
\ben
\mathbb{V}_h^k = \{ v\in L^2(\Omega): v|_K \in P^k(K), \forall K \in \mathcal{T}_h\}.
\een
Suppose $K$ and $K'$ are two adjacent triangles and share one common edge $e$. There are two traces of $v$ along the edge $e$, where we add or subtract those values to obtain the average and the jump.
The outward normal vector from $K$ to its neighbor element $K'$ is denoted by $\mathbf{n}=(n_1,n_2)$. Now the average and jump of $v$ on the edge $e$ are defined as follows,
$$
\overline {v}=\frac{1}{2}\left (v|_{K}+v|_{K'}\right ),\quad\quad [v]=v|_{K'}-v|_{K}.
$$

The original DDG scheme of (\ref{eqn:2DDiff}) defined in \cite{liu2009direct} is to find DG solution $u \in \mathbb{V}_h^k$, such that for any test function $v \in \mathbb{V}_h^k$ we have,
\begin{equation}\label{DDG_diffu_2D}
\int_K u_t v  ~dxdy  -\int_{\partial K} \widehat{\left(A(u)\nabla u \cdot {\bf n}\right)} v ~ds +
\int_K A(u) \nabla u \cdot \nabla v ~dxdy =0, \quad\forall \, K \in \mathcal{T}_h.
\end{equation}
The numerical flux $\widehat{A(u)\nabla u \cdot {\bf n}}$ (equation (3.7) of \cite{yan2013new} in dimension-by-dimension format) along the element edge is defined as,
$$
\widehat{A(u)\nabla u \cdot {\bf n}}
=\left(\widehat{b_{11}(u)_x}+\widehat{b_{12}(u)_y}\right)n_1
+\left(\widehat{b_{21}(u)_x}+\widehat{b_{22}(u)_y}\right)n_2,
$$
where $b_{ij}(u)=\int a_{ij}(u)du$ with $a_{ij}(u)$ as the diffusion matrix $A(u)$ entry. The outward normal is given with ${\bf n}=(n_1,n_2)$. Similar to (3.7) of \cite{yan2013new}, for example, the numerical flux $\widehat{b_{11}(u)_x}$ is calculated with formula,
$$
\widehat{b_{11}(u)_x}=\beta_{0} \frac{[b_{11}(u)]}{h_K}n_1 + \overline{b_{11}(u)_{x}} + \beta_1 h_K\left\{[b_{11}(u)_{xx}]n_1+[b_{11}(u)_{yx}]n_2\right\}.
$$
The coefficient pair $(\beta_0,\beta_1)$ should be chosen carefully to ensure the stability and convergence of the scheme. Again we have $h_K$ as element $K$'s diameter or the length of edge $\partial K$. Notice in \cite{Yan-MPS}, we further require the entries $a_{ij}(u)\geq 0$ to carry out the discrete maximum principle proof.

In this article, we simplify the calculation of numerical flux $\widehat{A(u)\nabla u \cdot {\bf n}}$ to the following,
\be\label{A_flux}
\widehat{A(u)\nabla u \cdot {\bf n}}  = \widehat{\nabla u \cdot \boldsymbol{\gamma}}
=\beta_0\frac{[u]}{h_K} + \overline{u_{\boldsymbol{\gamma}}} + \beta_1 h_K [u_{\boldsymbol{\gamma}\boldsymbol{\gamma}}],
\ee
where $\boldsymbol{\gamma} = A^T(u){\bf n}$ is a vector pointing from element $K$ into its neighbor along the edge $\partial K$. The simplification holds true since $(A(u)\nabla u) \cdot {\bf n}=\nabla u \cdot(A^T(u) {\bf n})$ and $A(u)$ is positive definite. We should point out the simplified version $(\ref{A_flux})$ is very important to carry out the maximum principle discussion in the following sections. Now the DDG interface correction \cite{liu2010direct} of (\ref{eqn:2DDiff}) is defined to find solution $u \in \mathbb{V}_h^k$, such that for any test function $v \in \mathbb{V}_h^k$ we have,
\be\label{DDG_IC_diffu_2D}
\int_K u_t v  ~dxdy  -\int_{\partial K} \widehat{\left(A(u)\nabla u \cdot {\bf n}\right)} v ~ds +
\int_K A(u) \nabla u \cdot \nabla v ~dxdy +
\int_{\partial K} \overline{A(v)\nabla v} \cdot {\bf n} [u] ~ds=0, \forall \, K \in \mathcal{T}_h,
\ee
with the numerical flux (\ref{A_flux}). The last term in (\ref{DDG_IC_diffu_2D}) is the
the extra added interface correction. Notice that the test function $v$ is taken to be zero outside the element $K$, thus the derivative average degenerates to $\overline{A(v)\nabla v}=\frac{1}{2}A(v)\nabla v|_K$ on the edge $\partial K$.

The complete scheme formulation of (\ref{eqn:2DConvecDiff}) with DDG interface correction follows,
\begin{align}\label{scheme:conv-diff}
\int_K u_tv ~dxdy+ &\int_{\partial K} \widehat{F\cdot{\bf n}} v ~ds - \int_K F\cdot\nabla v ~dxdy \notag\\
& =\int_{\partial K} \widehat{\left(A(u)\nabla u \cdot {\bf n}\right)} v ~ds -
\int_K A(u) \nabla u \cdot \nabla v ~dxdy -
\int_{\partial K} \overline{A(v)\nabla v} \cdot {\bf n} [u] ~ds,
\end{align}
with the convection term Lax-Friedrichs flux defined as,
\ben
\widehat{F\cdot{\bf n}} = \frac{1}{2}\left( F(u_K) \cdot {\bf n} + F(u_{K'}) \cdot{\bf n} - \alpha(u_{K'} - u_K) \right), \quad
\mbox{with}\quad \alpha = \max_{u} |{\bf F}'(u) \cdot {\bf n}|.
\een


\section{Maximum-Principle-Satisfying DDG methods}\label{sec:3rdMPSDDG}

In this section we prove  DDG polynomial solutions of nonlinear diffusion equations (\ref{eqn:2DDiff}) satisfy discrete maximum principle with the M-P-S limiter applied. We first discuss the linear case on unstructured triangle mesh in section \ref{sec:Linear diffusion equation}. Then we extend the study to nonlinear diffusion equations in section \ref{sec:Nonlinear diffusion equation}.

Notice that the second derivative jump term  has no contribution to the calculation of the numerical flux with low order $P^0$ and $P^1$ approximations. The scheme of DDG method with interface correction (\ref{DDG_IC_diffu_2D}) degenerates to IPDG methods with low order approximations. In this paper, we focus on $P^2$ quadratic polynomial approximations with 3rd order of accuracy. We skip the trivial piecewise constant case and refer to \cite{zhang-Shu-convdiff} for 2rd order linear approximations.

On the continuous level the maximum principle states that $m \leq u(x,y,t) \leq M$, given $m$ and $M$ as the minimum and maximum of the initial data $u(x,y,0)=u_0(x,y)$ and boundary data.
We have $u(x,y,t_n)$ to denote the exact solution at time level $t_n$ and $u^n_K(x,y)$ as the polynomial solution on element $K$ and at time level $t_n$. Our goal is to prove the polynomial solution satisfy $m\leq u^n_K(x,y)\leq M$ without losing the 3rd order accuracy at all time levels. We can simplify the discussion to Euler forward time discretization, since the full scheme (high order strong stability preserving (SSP) Runge-Kutta method) is a convex combination of Euler forward scheme.
For example, the third order SSP Runge-Kutta method in \cite{Shu-Osher-TVD-1988} is
\begin{equation*}
\left\{\begin{array}{lll}\label{3.31}
 u^{(1)}  = u^n+\Delta t H(u^n)\\
 u^{(2)}  = \frac{3}{4}u^n+\frac{1}{4}(u^{(1)}+\Delta t H(u^{(1)}))\\
u^{n+1}  = \frac{1}{3}u^n+\frac{2}{3}(u^{(2)}+\Delta t H(u^{(2)}))
\end{array}
\right.
\end{equation*}

\vspace{.05in}

Now assume at time level $t_n$, we have 1) the DDG solution is 3rd order accurate; and 2) we have $u_K^n(x,y)\in [m,M]$ on all elements. The goal is to prove the solution polynomial $u_K^{n+1}(x,y)$ at next time level $t_{n+1}$ still stay inside the bounds $[m,M]$ without losing accuracy. To carry out this study we need to consider following two steps:

\begin{enumerate}
\item
{\textit{to prove the polynomial solution average $\overline{u}_K^{n+1}$ stay inside the bounds $[m,M]$;}}
\item
{\textit{to prove the whole polynomial $u_K^{n+1}(x,y)$ stay inside $[m,M]$ without losing accuracy.}}
\end{enumerate}

The most challenging and the major step  is to show the polynomial average $\overline{u}_K^{n+1}$ falling in $[m,M]$. For the second step, we simply apply a linear scaling limiter \cite{Liu-Osher1996} to $u_K^{n+1}(x,y)$ and obtain a modified polynomial $\widetilde{u}_K^{n+1}(x,y)$ such that the whole polynomial $\widetilde{u}_K^{n+1}(x,y) \in [m,M]$ without losing accuracy. we refer to \cite{zhang-MPS-Review} for the proof of this accuracy preserving limiter. We have the DDG solution approximates the exact solution with 3rd order accuracy, thus the polynomial solution can only jump out of the bounds $[m,M]$ in the scale of $h^3$ with $h$ as the mesh size. The limiter \cite{Liu-Osher1996} is applied to compress and squeeze the polynomial in the scale of $h^{3}$ and put it back into the bounds $[m,M]$. The extra cost to preserve maximum principle is to apply the limiter \cite{Liu-Osher1996} to maintain the bounds.

\subsection{Linear diffusion equation}\label{sec:Linear diffusion equation}

In this section, we prove the DDG quadratic polynomial solutions of heat equation (\ref{eqn:2dheat}) satisfy discrete maximum principle. Again we focus on the first step and investigate under what conditions the solution average $\overline{u}_K^{n+1} \in [m,M]$ given $u_K^{n}(x,y)\in [m,M]$ on all elements.

The scheme formulation of DDG with interface correction for (\ref{eqn:2dheat}) is to find DG solution $u \in \mathbb{V}_h^k$, such that for any test function $v \in \mathbb{V}_h^k$ we have,
\begin{equation}\label{DDG_heat}
\int_K u_t v~dxdy  - \int_{\partial K} \widehat{u_{\textbf{n}}} v~ds+ \int_K \nabla u \cdot \nabla v ~dxdy +\int_{\partial K} \overline{v_{\textbf{n}}} [u]~ds= 0, \quad\forall \, K \in \mathcal{T}_h.
\end{equation}
The numerical flux $\widehat{u_{\bf n}}$ on the element boundary $\partial K$ is given with,
\be\label{DDG_heat_flux}
\widehat{u_{\bf n}} = \beta_0 \frac{[u]}{h_K} + \overline{u_{\bf n}} + \beta_1 h_K [u_{{\bf nn}}].
\ee
To obtain the solution average evolution, we take test function $v=1$ in (\ref{DDG_heat}), discretize in time with forward Euler and formally we have,
\begin{equation}{\label{avenext}}
\overline{u}^{n+1}_K = \overline{u}^n_K + \frac{\Delta t}{area(K)} \int_{\partial K}  \widehat{u_{\bf n}} ~ds.
\end{equation}
The solution average is given with $\overline{u}^{n}_K=\frac{1}{area(K)}\int_K u^n_K(x,y)~dxdy$ and $\Delta t$ is denoted as the time step size. The average $\overline{u}^{n}_K$ can be calculated out exactly, since $u^n_K(x,y)$ is a quadratic polynomial. Thus we see the quantity $\overline{u}^{n+1}_K$ of (\ref{avenext}) is essentially determined by the integral of $\widehat{u_{\bf n}}$ on the three edges. Recall that the numerical flux $\widehat{u_{\bf n}}$ of (\ref{DDG_heat_flux}) involves the solution jump, normal derivative average and second order normal derivative jumps on $\partial K$. The quantity $\overline{u}^{n+1}_K$ eventually is a function of the four solution polynomials that spread out in $K$ and its three neighbors.

For uniform triangular mesh in \cite{Yan-MPS}, we pick six solution values on each element and write out the $P^2$ solution polynomial in Lagrange format, then use a monotone argument to bound  $\overline{u}^{n+1}_K\in[m,M]$. For arbitrary triangular mesh, it's hard to identify such six points inside each triangle. We will use a new idea to bound $\overline{u}^{n+1}_K$.


We observe that the three quantities, namely $u$, $u_{\textbf{n}}$ and $u_{\textbf{nn}}$ restricted on the edges from each element, contribute to the calculation of $\widehat{u_{\textbf{n}}}$. If we manage to calculate $u$, $u_{\textbf{n}}$ and $u_{\textbf{nn}}$ from the given solution polynomials, we can easily bound $\overline{u}_K^{n+1}\in[m,M]$ once we have $u^n_K(x,y) \in[m,M]$ for all $K$.

\begin{figure}[htbp]
\centering
\includegraphics[width=0.45\linewidth]{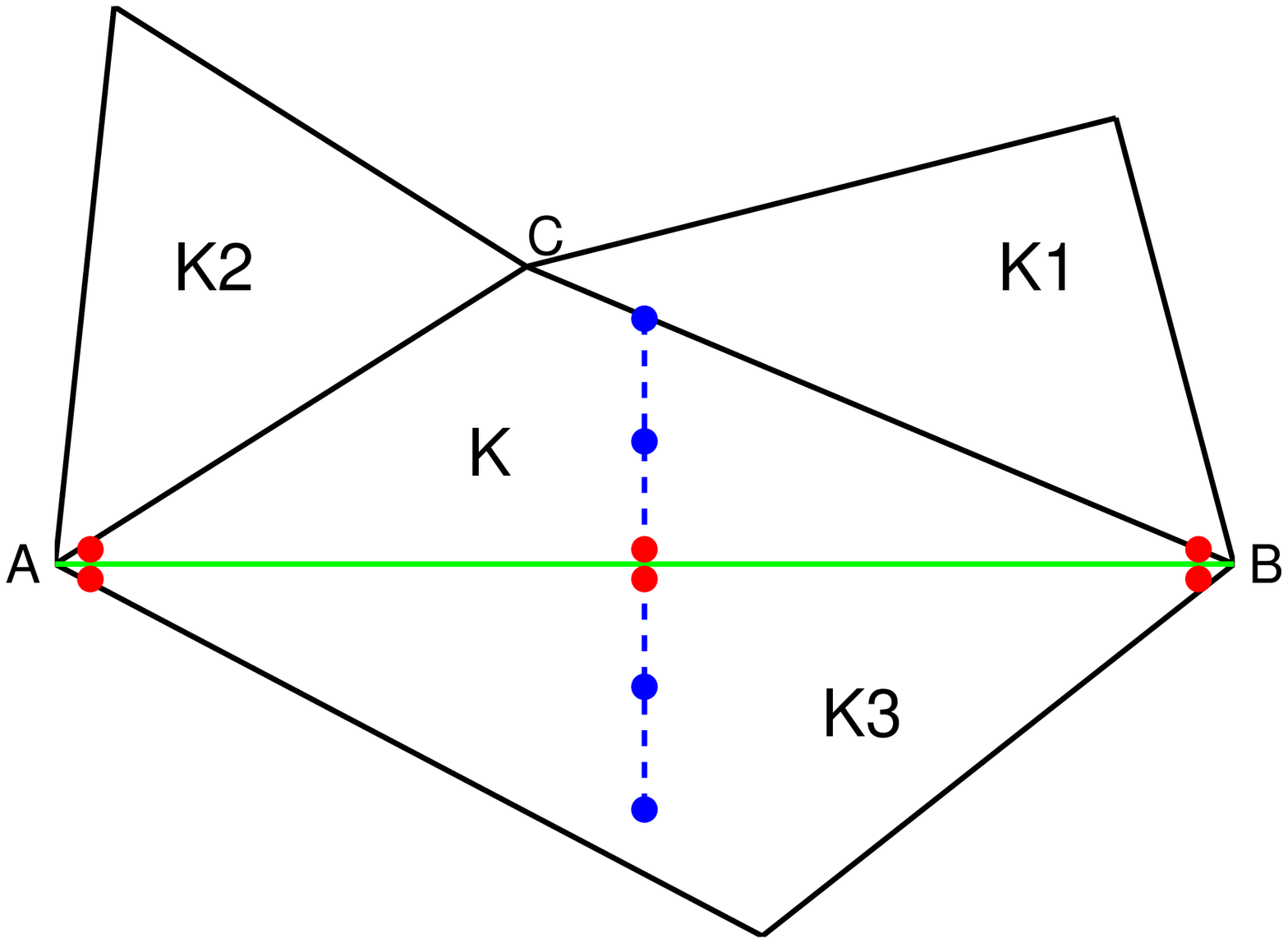}
\includegraphics[width=0.5\linewidth]{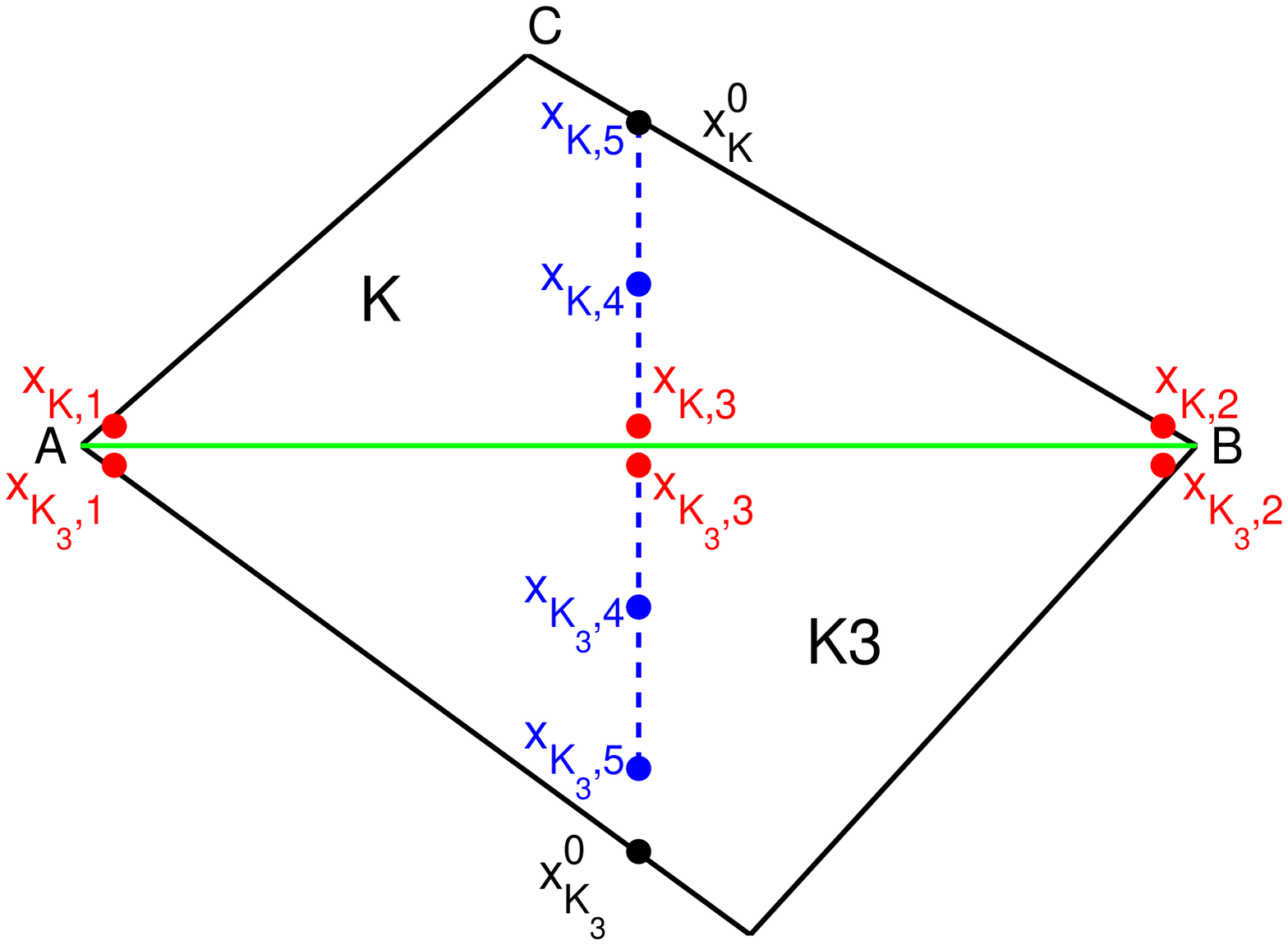}
\caption{Left: $K$ and its neighbor elements. Right: selected points to calculate $\widehat{u_{\bf n}}$ on edge $AB$.}
\label{fig:selectedpoints}
\end{figure}



Let's use Figure \ref{fig:selectedpoints} to illustrate the solution points selected to calculate $u$, $u_{\textbf{n}}$ and $u_{\textbf{nn}}$. For example, we consider the edge $AB$ shared by $K$ and its neighbor element $K_3$ (the right one in Figure \ref{fig:selectedpoints}). Notice that $u^n_{K_3}(x,y)$ is a $P^2$ polynomial. The three points, namely ${\bf x}_{K_3,1}$, ${\bf x}_{K_3,2}$, ${\bf x}_{K_3,3}$, are enough to represent the restriction of $u^n_{K_3}(x,y)$ on edge $AB$. The second normal derivative
$u_{\bf nn}$ degenerates to a constant, and we can use three points ${\bf x}_{K_3,3}^n$, ${\bf x}_{K_3,4}^n$, ${\bf x}_{K_3,5}^n$ on the normal line through the edge center to calculate $u_{\bf nn}$. The first normal derivative $u_{\bf n}$ on $AB$ is a linear polynomial, thus the same three points on the normal line are good enough to calculate the line integral of $u_{\bf n}$. With this new idea to calculate the numerical flux $\widehat{u_{\bf n}}$ of (\ref{avenext}), we are ready to bound the average $\overline{u}^{n+1}_K$.

\begin{thm}\label{thm:linear}
Consider DDG scheme with interface correction (\ref{DDG_heat}) - (\ref{DDG_heat_flux}) with $P^2$ quadratic approximations on unstructured triangular mesh. Given $u^n_K (x,y)$ in the range of $[m,M]$ for all $K$, we have $\overline{u}^{n+1}_K \in [m,M]$ provided,
\begin{equation}\label{cond}
 \beta_0 \geq \frac{9}{4} -6\beta_1, \quad \frac{1}{8} \leq \beta_1 \leq \frac{1}{4}, \quad \lambda= \frac{\Delta t}{aread(K)} \leq A(\beta_0, \beta_1, \check{\theta}, \hat{\theta}).
\end{equation}
Here $(\beta_0,\beta_1)$ is the coefficient pair in the numerical flux (\ref{DDG_heat_flux}).
We have $\check{\theta}$ and $\hat{\theta}$ denoted as the maximum and minimum angle of the partition $\mathcal{T}_h$, and $A$ is a function of $\beta_0$, $\beta_1$, $\check{\theta}$ and $\hat{\theta}$, i.e.,
\begin{equation}\label{A_cfl}
A=\tan(\hat{\theta})\cdot\min\left\{\frac{w_1}{72(1-4\beta_1)}, \frac{\tan(\hat{\theta})}{\tan(\check{\theta})}
\min\left(\frac{1}{6(8\beta_1-1)},\frac{w_1}{8(\beta_0-\frac{9}{4}+6\beta_1)},\frac{w_1}{4\beta_0}\right)
\right\},
\end{equation}
where $w_1=\frac{2}{81}$ as shown in (\ref{eq:weight4vertex}).
\end{thm}

\vspace{.1in}

\begin{proof}

To bound  $\overline{u}^{n+1}_K$ of (\ref{avenext}), we see it's important to carefully calculate the numerical flux $ \widehat{u_{\bf n}}$ on $\partial K$. From the numerical flux formula (\ref{DDG_heat_flux}) of the DDG schemes, we have $h_K$ taken as the element diameter or the length of the edge $\partial K$. To simplify the proof, here we modify $h_K$ to incorporate with the mesh geometrical information. We should comment that numerically we observe no difference with either choice of $h_K$.

Again we use edge $AB$ in Figure \ref{fig:selectedpoints} to illustrate the definition of $h_K$ chosen in the numerical flux formula (\ref{DDG_heat_flux}). Let's have the parametric equation ${\bf r}(t)=t{\bf n}+{\bf x}_{K,3}$, $t\in\mathbb{R}$ to represent the normal line through the edge center. And we have points ${\bf x}_{K}^0$ and ${\bf x}_{K_3}^0$ as the intersection of the normal line with the other two edges of $K$ and $K_3$. Restricted on edge $\partial K=AB$, we take $h_K=h_{_{AB}} = \min\left\{\|{\bf x}_{K,3}-{\bf x}_{K}^0\|,\|{\bf x}_{K_3,3}-{\bf x}_{K_3}^0\|\right\}$ and we have,

\begin{equation}\label{h_{AB}}
\int_{AB} \widehat{u_{\bf n}} ~ds = \int_{AB}  \beta_0 \frac{[u]}{h_{_{AB}}} ~ds + \int_{AB}  \overline{u_{\bf n}} ~ds + \int_{AB} \beta_1 h_{_{AB}} [u_{{\bf nn}}] ~ds.
\end{equation}


To calculate $u_{{\bf nn}}$, we pick two more points along the normal line from each side. We have points ${\bf x}_{K,5} (t=-h_{_{AB}})$ and ${\bf x}_{K,4} (t=-\frac{1}{2}h_{_{AB}})$ taken in element $K$, and  points ${\bf x}_{K_3,5} (t=h_{_{AB}})$ and ${\bf x}_{K_3,4} (t=\frac{1}{2}h_{_{AB}})$ taken in element $K_3$, as shown in Figure \ref{fig:selectedpoints}. Furthermore we denote $u_{K,1},\dots, u_{K,5}$ as the $u^n_K(x,y)$ quadratic polynomial solution values on points ${\bf x}_{K,1},\dots, {\bf x}_{K,5}$. As discussed previously, the five points $u_{K,1},\dots, u_{K,5}$ are enough to calculate $\int_{AB} \widehat{u_{\bf n}} ds$ from the side of element $K$. Similarly, the five points $u_{K_3,1},\dots, u_{K_3,5}$ are enough to calculate $\int_{AB} \widehat{u_{\bf n}} ds$ from the side of element $K_3$. Essentially the quantity  $\int_{AB} \widehat{u_{\bf n}} ds$ can be explicitly written out in terms of the ten solution values spread out in elements $K$ and $K_3$ as follows,

\begin{align}\label{un_edge_AB}
\int_{AB} \beta_0 \frac{[u]}{h_{_{AB}}} ds
&=  \frac{\beta_0l_{_{AB}}}{6h_{_{AB}}}
\left\{\left(u_{K_3,1}+u_{K_3,2}+4u_{K_3,3}\right) - \left(u_{K,1}+u_{K,2}+4u_{K,3} \right) \right\}\notag\\
\int_{AB} \overline{u_{\bf n}} ds
&= \frac{l_{_{AB}}}{2h_{_{AB}}}\left\{\left(-3u_{K_3,3} - u_{K_3,5}+4u_{K_3,4}\right) + \left(3u_{K,3} + u_{K,5}-4u_{K,4}\right)\right\} \notag\\
\int_{AB} \beta_1 h_{_{AB}} [u_{{\bf n}{\bf n}}] ds
&=  \frac{4\beta_1 l_{_{AB}}}{h_{_{AB}}}\left\{\left(u_{K_3,3} + u_{K_3,5}-2u_{K_3,4}\right) - \left(u_{K,3} + u_{K,5}-2u_{K,4}\right)\right\}.\notag\\
\end{align}
Here we have $l_{_{AB}}$ denoting the length of edge $AB$. Finally the average $\overline{u}^{n+1}_K$ of (\ref{avenext}) can be written out as a function of the solution values that spread out in element
$K$ and its neighbors $K_1, K_2, K_3$. Thus formally we have,
\begin{align}\label{eq:ave_Hfunc}
\overline{u}^{n+1}_K
&= \overline{u}^n_K + \frac{\Delta t}{area(K)} \left\{ \int_{AB} \widehat{u_{\bf n}} ds + \int_{BC}\widehat{u_{\bf n}} ds + \int_{CA} \widehat{u_{\bf n}}ds \right\} \notag\\
&= H\left\{u^n_K(\cdot,\cdot),u^n_{K_1}(\cdot,\cdot),u^n_{K_2}(\cdot,\cdot),u^n_{K_3}(\cdot,\cdot)\right\}.
\end{align}
The functional $H(\cdots)$ involves 28 arguments with 15 points from $K_1, K_2, K_3$ and 13 points from $K$. The first 12 points from $K$ are selected from the calculation of $\int_{\partial K} \widehat{u_{\bf n}} ds$, see Figure \ref{fig:selectedpoints}. The 13th one is to be selected by the quadrature rule for cell average $\overline{u}^n_K$, see Appendix \ref{quadrature}.

Our goal here is to prove solution average $\overline{u}^{n+1}_K\in[m,M]$ given $u^n_K (x,y)\in [m,M]$ on all elements. Again, we use a monotone argument showing  $\overline{u}^{n+1}_K$ is a convex combination of the selected solution points. To study the conditions to guarantee that $H(\Uparrow,\Uparrow,\Uparrow,\Uparrow)$ is monotonically increasing on the total 28 arguments, it is enough to check out the ten points selected inside $K$ and $K_3$ of (\ref{eq:ave_Hfunc}). We first check the five points selected on element $K_3$. From (\ref{un_edge_AB}) - (\ref{eq:ave_Hfunc}) we have,
$$
\frac{\partial H}{\partial u_{K_3,1}} = \frac{\partial H}{\partial u_{K_3,2}}
= \lambda\frac{l_{_{AB}}}{h_{_{AB}}}\frac{\beta_0}{6}, \quad\quad
\frac{\partial H}{\partial u_{K_3,3}}
= \lambda\frac{l_{_{AB}}}{h_{_{AB}}}(\frac{2}{3}\beta_0 - \frac{3}{2} + 4\beta_1),
$$

$$
\frac{\partial H}{\partial u_{K_3,4}} = \lambda\frac{2l_{_{AB}}}{h_{_{AB}}}(1 - 4\beta_1), \quad\quad
\frac{\partial H}{\partial u_{K_3,5}}= \lambda\frac{l_{_{AB}}}{2h_{_{AB}}}(8\beta_1-1).
$$
With $\lambda = \frac{\Delta t}{area(K)}>0$, we only need $ \beta_0 \geq \frac{9}{4} -6\beta_1$ and
$\frac{1}{8} \leq \beta_1 \leq \frac{1}{4}$ to guarantee the coefficients of the five solution values in $K_3$  being non-negative. Before we carry out the discussion on the five points in element $K$, we need an inequality (refer to Appendix \ref{quadrature} ) which reflects geometrical property of the mesh partition.
\begin{equation*}
h_{_{AB}}=\frac{l_{_{AB}}}{2}\tan(\min(\theta_1,\theta_2,\theta_4,\theta_5)) \geq \frac{l_{_{AB}}}{2}\tan(\hat{\theta}).
\end{equation*}
From (\ref{un_edge_AB}) and (\ref{eq:ave_Hfunc}), we have,
\begin{align*}
\frac{\partial H}{\partial u_{K,1}}
&=\frac{\partial \overline{u}^n_K}{\partial u_{K,1}} -  \lambda\frac{\beta_0}{6}  \left(\frac{l_{_{AB}}}{h_{_{AB}}} +\frac{l_{_{CA}}}{h_{_{CA}}} \right)
\geq  \frac{\partial \overline{u}^n_K}{\partial u_{K,1}} -  \lambda\frac{2\beta_0}{3 \tan(\hat{\theta})}\notag\\
\frac{\partial H}{\partial u_{K,2}}
&= \frac{\partial \overline{u}^n_K}{\partial u_{K,2}} -  \lambda\frac{\beta_0}{6}
\left(\frac{l_{_{AB}}}{h_{_{AB}}} +\frac{l_{_{BC}}}{h_{_{BC}}}\right)
\geq  \frac{\partial \overline{u}^n_K}{\partial u_{K,2}} -  \lambda\frac{2\beta_0}{3 \tan(\hat{\theta})}\notag\\
\frac{\partial H}{\partial u_{K,3}}
&= \frac{\partial \overline{u}^n_K}{\partial u_{K,3}} -  \lambda\frac{l_{_{AB}}}{h_{_{AB}}}\left(\frac{2}{3}\beta_0 - \frac{3}{2} + 4\beta_1\right)
\geq \frac{\partial \overline{u}^n_K}{\partial u_{K,3}} -  \lambda\frac{4\left(\beta_0 - \frac{9}{4} + 6\beta_1\right)}{3\tan(\hat{\theta})}\notag\\
\frac{\partial H}{\partial u_{K,4}}
&= \frac{\partial \overline{u}^n_K}{\partial u_{K,4}} - \lambda\frac{2l_{_{AB}}}{h_{_{AB}}}\left(1 - 4\beta_1\right)
\geq  \frac{\partial \overline{u}^n_K}{\partial u_{K,4}} - \lambda\frac{4\left(1 - 4\beta_1\right)}{\tan(\hat{\theta})}\notag\\
\frac{\partial H}{\partial u_{K,5}}
&= \frac{\partial \overline{u}^n_K}{\partial u_{K,5}} - \lambda\frac{l_{_{AB}}}{2h_{_{AB}}}\left(8\beta_1-1\right)
\geq   \frac{\partial \overline{u}^n_K}{\partial u_{K,5}} - \lambda\frac{ 8\beta_1-1}{\tan(\hat{\theta})}.  \notag\\
\end{align*}
To guarantee the coefficients of the five solution points in $K$ being non-negative, we need a special quadrature rule with all positive weights on all the 12 selected points in $K$. We refer to Appendix \ref{quadrature} for details of this quadrature rule. The 13th point $u_{K,13}$ is selected by the quadrature rule also with positive weight. With CFL restriction (\ref{cond})- (\ref{A_cfl}) and the quadrature weights (\ref{weights}), we see $H(\cdots)$ is monotonically increasing on $u_{K,1}$, $u_{K,2}$, $u_{K,3}$, $u_{K,4}$ and $u_{K,5}$.
Similar argument applies to edge $BC$ and edge $CA$, involving solution values in elements $K_1$ and $K_2$. Easily we see functional $H(\cdots)$ is monotonically increasing w.r.t. all 28 point values.
With the consistency and the monotonicity of $H(\cdots)$, we obtain,
$$
m=H(m,\cdots,m)\leq \overline{u}_K^{n+1} =H(\cdots)\leq H(M,\cdots,M)=M.
$$
provided that $m\leq u^n_K(x,y), u^n_{K_1}(x,y), u^n_{K_2}(x,y), u^n_{K_3}(x,y)\leq M$.


\end{proof}

\subsection{Nonlinear diffusion equation}\label{sec:Nonlinear diffusion equation}

In this section we extend the study of 3rd order M-P-S DDG scheme with interface correction (\ref{A_flux}) - (\ref{DDG_IC_diffu_2D}) to general nonlinear diffusion equations (\ref{eqn:2DDiff}) on unstructured triangle mesh. Again, our goal is to bound the solution average $\overline{u}^{n+1}_K \in [m,M]$ provided $u^n_K (x,y)$ in the range of $[m,M]$. Take test function $v=1$ in (\ref{DDG_IC_diffu_2D}) and discretize in time with Euler forward, we have the solution average evolving in time as,
\be\label{2DAverageEvo-2}
\overline{u}^{n+1}_K = \overline{u}^n_K + \frac{\Delta t}{area(K)}\int_{\partial K} \widehat{\left(A(u)\nabla u \cdot {\bf n}\right)}  ~ds=\overline{u}^n_K + \frac{\Delta t}{area(K)}\int_{\partial K} \widehat{\nabla u \cdot \boldsymbol{\gamma} }~ds.
\ee

\begin{figure}[htbp]
\centering
\includegraphics[width=0.5\linewidth]{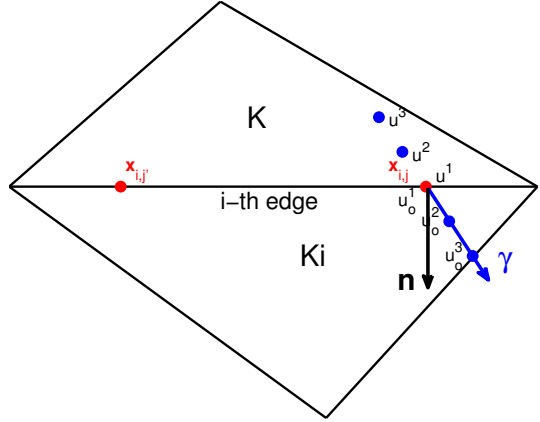}
\caption{Selected points along direction $\gamma = {\bf A}^T(u^n({\bf x}_{i,j})){\bf n}_i$ for representing the numerical flux}
\label{fig:nonlinear_gauss_points}
\end{figure}
As discussed previously in (\ref{A_flux}), we apply a new way to calculate numerical flux $\widehat{A(u)\nabla u \cdot {\bf n}}=\widehat{\nabla u \cdot \boldsymbol{\gamma} }$ with $\boldsymbol{\gamma} = A^T(u){\bf n}$ as a vector pointing from $K$ into its neighbor along the edge. This is true because the diffusion matrix $A(u)$ is positive definite and we have
$\boldsymbol{\gamma} \cdot  {\bf n}= A^T(u){\bf n} \cdot {\bf n} > 0$. In a word, $\gamma$ is a vector always pointing into its neighbor element. To bound $\overline{u}^{n+1}_K$, we need to manage to calculate $\int_{\partial K} \widehat{\nabla u \cdot \boldsymbol{\gamma} }~ds$ on the three edges of element $K$. Notice that vector $\boldsymbol{\gamma} $ is a nonlinear function of the solution, thus we apply a quadrature rule to calculate $\int_{\partial K} \widehat{\nabla u \cdot \boldsymbol{\gamma} }~ds$. For example, we consider 2-point Gaussian quadrature rule along each edge to approximate the line integral. Let's use point ${\bf x}_{i,j}$ to denote the $j$th Gaussian point on the $i$th edge. For each Gaussian point ${\bf x}_{i,j}$, shown in Figure \ref{fig:nonlinear_gauss_points}, we see six solution values are enough to calculate the numerical flux,
\be\label{flux_on_gauss_point}
\widehat{\nabla u \cdot \boldsymbol{\gamma}}|_{{\bf x}_{i,j}}=\widehat{u_{\boldsymbol{\gamma}}}|_{{\bf x}_{i,j}}
=  \beta_0\frac{[u]}{h_{i,j}} + \overline{u_{\boldsymbol{\gamma}}} + \beta_1 h_{i,j} [u_{\boldsymbol{\gamma\gamma}}]\biggr\rvert_{{\bf x}_{i,j}},
\ee
where $\boldsymbol{\gamma}= A^T(u^n({\bf x}_{i,j})){\bf n}_i$ and $h_{i,j}$ as the shortest distance from point ${\bf x}_{i,j}$ to the other edges of $K$ and $K_i$ along $\boldsymbol{\gamma}$. Again, we bound $\overline{u}^{n+1}_K$ by showing it is a convex combination of polynomial solution values that spread in element $K$ and its three neighbors $K_1$, $K_2$ and $K_3$.


\begin{thm}(Nonlinear diffusion equation)
Consider DDG scheme with interface correction (\ref{A_flux}) - (\ref{DDG_IC_diffu_2D}) with $P^2$ quadratic approximations on a triangular mesh. Given $u^n_K (x,y)$ in the range of $[m,M]$, we have $\overline{u}^{n+1}_K \in [m,M]$ provided,
\begin{equation}\label{cond2}
 \beta_0 \geq \frac{3}{2} -4\beta_1, \quad \frac{1}{8} \leq \beta_1 \leq \frac{1}{4}, \quad \lambda= \frac{\Delta t}{area(K)} \leq A(\beta_0, \beta_1, \check{\theta}, \hat{\theta}).
\end{equation}
Again $(\beta_0,\beta_1)$ is the coefficient pair in the numerical flux (\ref{A_flux}), and $\check{\theta}$ and $\hat{\theta}$ as the maximum and minimum angles of the partition $\mathcal{T}_h$. Function $A$ depends on $\beta_0$, $\beta_1$, $\check{\theta}$ and $\hat{\theta}$ as,
\begin{equation}\label{A_cfl_nonlinear}
A=\sin(\hat{\theta})\cdot \frac{3-\sqrt{3}}{3}w_0\cdot\min\left\{
\frac{1}{2\beta_0+8\beta_1+3}, \frac{1}{8\beta_1+1}
\right\},
\end{equation}
where $w_0$ is the minimum quadrature weight in the quadrature rule.
\end{thm}

\vspace{.05in}

\begin{proof}

Similar to Theorem \ref{thm:linear}, it suffices to study the monotonicity of $\overline{u}^{n+1}_K$ with respect to the points selected to evaluate the right hand side of (\ref{2DAverageEvo-2}). Specifically it is enough to study the Gaussian points that are used to approximate the line integral on the edges.
For one Gaussian point ${\bf x}_{i,j}$, as shown in Figure \ref{fig:nonlinear_gauss_points}, six solution points selected along $\gamma = A^T(u^n({\bf x}_{i,j})){\bf n}_i$ direction are enough to calculate the numerical flux (\ref{flux_on_gauss_point}). We denote  $u^1$, $u^2$, $u^3$ as the three solution values of $u^n_K(x,y)$ selected in $K$ and $u^1_o$, $u^2_o$, $u^3_o$ as the ones selected in $K_i$. Again we have $h_{i,j}$ denoting the shortest distance from point ${\bf x}_{i,j}$ to the other edges of $K$ and $K_i$ along $\boldsymbol{\gamma}$.

Now at the Gaussian point ${\bf x}_{i,j}$, we can explicitly write out the value of $\widehat{u_{\boldsymbol{\gamma}}}$ of (\ref{flux_on_gauss_point}) with,
\begin{align}
&[u] = u^1_o - u^1 \notag \\
&\overline{u_{\gamma}} = \frac{1}{2}\left(\frac{4u^2-3u^1-u^3}{h_{i,j}} + \frac{4u^2_o-3u^1_o-u^3_o}{h_{i,j}}\right)
=\frac{1}{2h_{i,j}}\left(4u^2-3u^1-u^3 + 4u^2_o-3u^1_o-u^3_o\right) \notag\\
&[u_{\gamma\gamma}]= \frac{u^1_o + u^3_o -2u^2_o}{\left(\frac{h_{i,j}}{2}\right)^2} - \frac{u^1 + u^3 -2u^2}{\left(\frac{h_{i,j}}{2}\right)^2}
= \frac{4}{h_{i,j}^2} \left[u^1_o + u^3_o -2u^2_o- (u^1 + u^3 -2u^2) \right].\notag
\end{align}
Introduce notation $\alpha^{i,j} =  \frac{\Delta t}{area(K)}l_K^i \omega_j$ with $l_K^i$ as the $i$th edge length and $\omega_j$ as the  $j$th Gaussian point weight. Since we use 2-point Gaussian quadrature rule, the weight $\omega_j \equiv 1$. From  (\ref{2DAverageEvo-2}), we have,
\begin{align*}
&\frac{\partial \overline{u}^{n+1}_K}{\partial u^3_o}  = \frac{\alpha^{i,j}}{2h_{i,j}} (8\beta_1 -1),\quad
\frac{\partial \overline{u}^{n+1}_K}{\partial u^2_o}  =  \frac{2\alpha^{i,j}}{h_{i,j}} (1 - 4\beta_1),\,\notag \\
&\frac{\partial \overline{u}^{n+1}_K}{\partial u^1_o}  = \frac{\alpha^{i,j}}{h_{i,j}} (\beta_0 +4\beta_1 - \frac{3}{2}),\quad
\frac{\partial \overline{u}^{n+1}_K}{\partial u^3}  = \frac{\partial \overline{u}^{n}_K}{\partial u^3} - \frac{\alpha^{i,j}}{2h_{i,j}}(8\beta_1 + 1), \,\notag \\
&\frac{\partial \overline{u}^{n+1}_K}{\partial u^2}  = \frac{\partial \overline{u}^{n}_K}{\partial u^2} + \frac{2\alpha^{i,j}}{h_{i,j}}( 4\beta_1 + 1),\,\quad
\frac{\partial \overline{u}^{n+1}_K}{\partial u^1} = \frac{\partial \overline{u}^{n}_K}{\partial u^1} - \frac{\alpha^{i,j}}{2h_{i,j}}( 2\beta_0 + 8\beta_1 + 3).
\end{align*}
With $\beta_0 \geq \frac{3}{2} -4\beta_1, \, \frac{1}{8} \leq \beta_1 \leq \frac{1}{4}$ satisfied in the numerical flux (\ref{A_flux}), we have  $\overline{u}^{n+1}_K$ as a monotone increasing function on the solution values $u^1_o$, $u^2_o$, $u^3_o$ chosen from element $K_3$. Similar discussion applies to the solution values used in element $K_1$ and $K_2$. With 2-point Gaussian quadrature rule approximating the line integral, the quantity  $\overline{u}^{n+1}_K$ is monotone increasing with respect to the total 18 solution values from $K_1$, $K_2$ and $K_3$.

Similar to Theorem \ref{thm:linear},  we need to have a special quadrature rule for the cell average $\overline{u}^n_K$ that use all the selected points inside $K$ (18 points in total) with positive weights. We use a similar method as the linear case to find such a quadrature rule, see Appendix \ref{quadrature}. Let $w_0$ be the minimum weight for the selected points.

With geometry information of the mesh partition,
\be
\frac{l_K^i}{h_{i,j}} \le \frac{1}{C\sin{\hat{\theta}}}, \qquad C \ge \frac{3-\sqrt{3}}{6},
\ee
we see condition (\ref{cond2}) is sufficient to guarantee the coefficients of solution values used in element $K$ being none-negative. Therefore, $\overline{u}^{n+1}_K$ is monotonically increasing w.r.t. all the selected points inside $K$ and its neighbors. Finally we conclude with,
$$
m\leq \overline{u}_K^{n+1}\leq M,
$$
provided that $m\leq u^n_K(x,y), u^n_{K_1}(x,y), u^n_{K_2}(x,y), u^n_{K_3}(x,y)\leq M$.

\end{proof}

\vspace{.2in}
\noindent
{\bf \large Implementation of the M-P-S limiter}

Given the quadratic DG polynomial solution $u^n_K(x,y)$ with cell average $\overline{u}_K^n \in [m, M]$, the following limiter ensures $\widetilde{u}^n_K(x,y) \in [m, M]$ for any $(x,y) \in K$.
\be\label{limiter}
\widetilde{u}^n_K(x,y) = \theta (u^n_K(x,y) - \overline{u}_K^n) + \overline{u}_K^n, \qquad
\theta = \min\left\{ 1, \left|\frac{M - \overline{u}_K^n}{M_K - \overline{u}_K^n}\right|, \left|\frac{m - \overline{u}_K^n}{m_K - \overline{u}_K^n}\right|	\right\},
\ee
with $M_K$ and $m_K$ as the maximum and minimum of $u^n_K(x,y)$ on element $K$,
\be \label{max_min}
M_K = \max_{(x,y) \in K} u^n_K(x,y), \qquad m_K = \min_{(x,y) \in K} u^n_K(x,y).
\ee
Since $u^n_K(x,y)$ is a quadratic polynomial, it is easy to calculate the maximum and minimum value over $K$.

Notice that the limiter (\ref{limiter}) does not change the cell average. Moreover, $\| \widetilde{u}^n_K(x,y) - u^n_K(x,y)\|_\infty = \mathcal{O}(h^3)$, if the exact solution is smooth. The proof can be found in \cite{zhang2010maximum}. Thus, $\widetilde{u}^n_K(x,y) \in [m, M]$ has uniform third order accuracy for smooth exact solution.

\begin{Algorithm} {\bf maximum-principle-satisfying DDG scheme with interface correction} \label{algorithm}
\noindent
\begin{enumerate}
\item
At time level $t_n$, we apply M-P-S limiter (\ref{limiter}) - (\ref{max_min}) to $u^n_K(x,y)$ and obtain $\widetilde{u}_K^{n}(x,y)$.
\item
Apply DDG scheme with interface correction  (\ref{A_flux}) - (\ref{DDG_IC_diffu_2D}) to $\widetilde{u}_K^{n}(x,y)$ and evolve in time with SSP Runge-Kutta method \cite{Shu-Osher-TVD-1988} to march forward the solution to the next time level $t_{n+1}$.
\end{enumerate}
\end{Algorithm}

For the convection part of (\ref{eqn:2DConvecDiff}), as discussed in \cite{zhang2010maximum}, the solution average at next time level $\overline{u}_K^{n+1}$ is a monotone function with respect to certain solution values (Gauss-Lobatoo points) at time level $t_n$. Thus the  M-P-S limiter (\ref{limiter}) - (\ref{max_min}) with $M_K$ and $m_K$ as the maximum and minimum of the polynomial solution $u^n_K(x,y)$ over the whole element $K$ is enough to guarantee the solution average staying in the given bound.

\begin{rem}
For general convection diffusion equation (\ref{eqn:2DConvecDiff}), we apply same procedure listed in Algorithm \ref{algorithm} to guarantee the quadratic polynomial solution stay in the given bound and at the same time maintain the 3rd order accuracy. For (\ref{eqn:2DConvecDiff}), we apply DDG with interface correction scheme (\ref{scheme:conv-diff}) for spatial discretization.
\end{rem}

\section{Numerical Examples}\label{sec:numerical}
\numberwithin{equation}{section}

In this section, we present a sequence of examples to demonstrate the performance of M-P-S limiter.
For all examples in this section, we have the coefficient pair taken as $(\beta_0,\beta_1)=(5,\frac{1}{8})$ in the numerical flux formula (\ref{DDG_heat_flux}).

\begin{example}\label{e.g.:accuracy_test} {Accuracy test on linear diffusion equation}
\end{example}

\begin{figure}[!htbp]
\centering
\subfigure[Triangular mesh with $h = 0.117$.]{\includegraphics[width=0.45\textwidth]{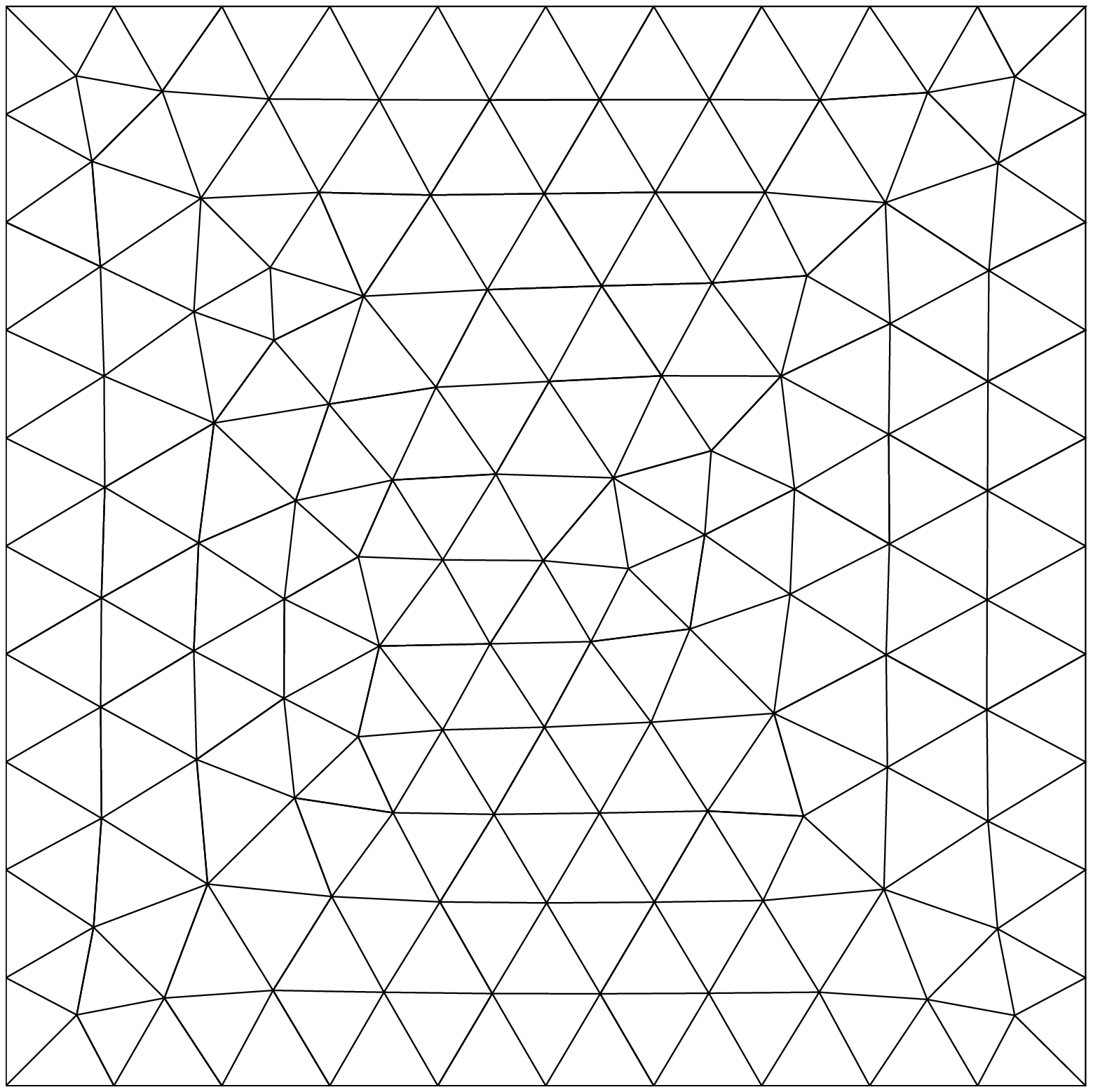}\label{fig:grid_3}}
\subfigure[Mesh with obtuse triangles, $h = 0.148$.]{\includegraphics[width=0.45\textwidth]{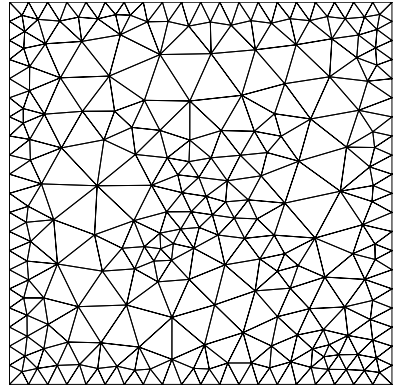}\label{fig:grid_5}}
\caption{Illustration of meshes}
\label{fig:quadS2K}
\end{figure}

We start with accuracy check of the DDG with interface correction (\ref{DDG_heat}) with and without M-P-S limiter (\ref{limiter}) applied on the following linear diffusion equation,
\ben
u_t  - \epsilon\Delta u = 0, \quad\quad  (x,y) \in [0,1]\times[0,1],\qquad t \in (0,T),
\een
with initial data $u(x,y,0) = u_0(x,y) = \sin(2\pi(x+y))$ and periodic boundary condition. The exact solution is given with,
\ben
u(x,y,t) = \exp(-8\pi^2\epsilon t)\sin(2\pi (x+y)).
\een
Here, we take $\epsilon = 1$ and final time $t=0.0001$.  We implement the scheme with $P^2$ quadratic polynomials on unstructured mesh in Figure \ref{fig:grid_3} and on mesh with obtuse triangles with largest angle about $\frac{3}{5}\pi$ in Figure \ref{fig:grid_5}. Third order of accuracy is maintained with and without M-P-S limiter (\ref{limiter}) applied, see Table \ref{tab:accuracy_3} and Table \ref{tab:accuracy_5}. At each time step $t_n$, we set the bounds to be $ue_{\min} = -\exp(-8\pi^2\epsilon t_n)$ and $ue_{\max} = \exp(-8\pi^2\epsilon t_n)$, which are the minimum and maximum of the exact solution. We use $u_{\min}$ and $u_{\max}$ to denote the DG solution minimum and maximum values. The overshoots and undershoots are eliminated after the M-P-S limiter applied, see Table \ref{tab:accuracy_3} and Table \ref{tab:accuracy_5}.

\begin{table}[htbp]
   \centering
\begin{minipage}{1.\textwidth} 
    \resizebox{\textwidth}{!}{
\begin{tabular}{|c|cccc|cc|cccc|cc|}
\hline
 $h$& \multicolumn{6}{c|}{without M-P-S limiter} & \multicolumn{6}{c|}{with M-P-S limiter} \\
\cline{2-13}
 & $L^2$ error & Order & $L^\infty$ error & Order & $u_{\min}-ue_{\min}$ & $u_{\max}-ue_{\max}$ & $L^2$ error & Order & $L^\infty$ error & Order & $u_{\min}-ue_{\min}$ & $u_{\max}-ue_{\max}$\\
 \hline
0.117	&	1.21e-03	&		&	1.10e-02	&		&	-3.91e-03	&	3.16e-03	&	1.28e-03	&		 &	1.10e-02	 &		 &	0	&	0	\\ \hline
0.0587	&	1.94e-04	&	2.64	&	1.28e-03	&	3.10	&	-2.40e-04	&	2.49e-04	&	 1.95e-04	&	2.71	&	 1.28e-03	&	3.10	&	0	&	0	\\ \hline
0.0293	&	2.60e-05	&	2.90	&	1.59e-04	&	3.01	&	-1.49e-05	&	1.50e-05	&	 2.60e-05	&	2.91	&	 1.59e-04	&	3.01	&	0	&	0	\\ \hline
0.0147	&	3.27e-06	&	2.99	&	2.01e-05	&	2.98	&	-9.61e-07	&	9.41e-07	&	 3.27e-06	&	2.99	&	 2.01e-05	&	2.98	&	0	&	0	\\ \hline
0.00733	&	4.11e-07	&	2.99	&	2.52e-06	&	2.99	&	-5.82e-08	&	5.96e-08	&	 4.11e-07	&	2.99	&	 2.52e-06	&	2.99	&	0	&	0	\\ \hline
\end{tabular}}
\caption{Accuracy table on triangular mesh of Figure \ref{fig:grid_3}.}
\label{tab:accuracy_3}
   \end{minipage}
\vspace{0.5cm}

   \begin{minipage}{1.\textwidth} 
    \resizebox{\textwidth}{!}{
\begin{tabular}{|c|cccc|cc|cccc|cc|}
\hline
 $h$& \multicolumn{6}{c|}{without M-P-S limiter} & \multicolumn{6}{c|}{with M-P-S limiter} \\
\cline{2-13}
 & $L^2$ error & Order & $L^\infty$ error & Order & $u_{\min}-ue_{\min}$ & $u_{\max}-ue_{\max}$ & $L^2$ error & Order & $L^\infty$ error & Order & $u_{\min}-ue_{\min}$ & $u_{\max}-ue_{\max}$\\
 \hline
0.148	&	8.23e-04	&		&	1.46e-02	&		&	-1.17e-02	&	9.70e-03	&	1.08e-03	&		 &	1.46e-02	 &		 &	0	&	0	\\ \hline
0.0741	&	1.23e-04	&	2.74	&	1.82e-03	&	3.00	&	-7.86e-04	&	5.49e-05	&	 1.26e-04	&	3.11	&	 2.02e-03	&	2.85	&	0	&	0	\\ \hline
0.0371	&	1.68e-05	&	2.87	&	2.17e-04	&	3.06	&	-2.66e-05	&	4.24e-05	&	 1.68e-05	&	2.90	&	 2.17e-04	&	3.22	&	0	&	0	\\ \hline
0.0185	&	2.14e-06	&	2.97	&	2.69e-05	&	3.02	&	-1.12e-06	&	9.22e-07	&	 2.14e-06	&	2.97	&	 2.69e-05	&	3.02	&	0	&	0	\\ \hline
0.00927	&	2.70e-07	&	2.99	&	3.33e-06	&	3.01	&	-7.17e-08	&	1.33e-08	&	 2.70e-07	&	2.99	&	 3.33e-06	&	3.01	&	0	&	0	\\ \hline
\end{tabular}}
\caption{Accuracy table on unstructured mesh with obtuse triangles of Figure \ref{fig:grid_5}.}
\label{tab:accuracy_5}
   \end{minipage}
\end{table}

\begin{example}\label{e.g.:porous} {Porous medium equation}
\end{example}

\begin{figure}
\centering
\subfigure[t = 0]{\includegraphics[width=0.45\textwidth]{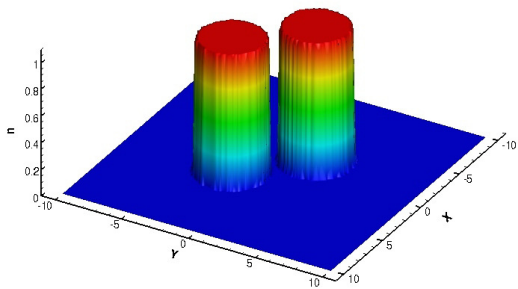}\label{fig:porous_5_t_00}}
\subfigure[t = 0.1]{\includegraphics[width=0.45\textwidth]{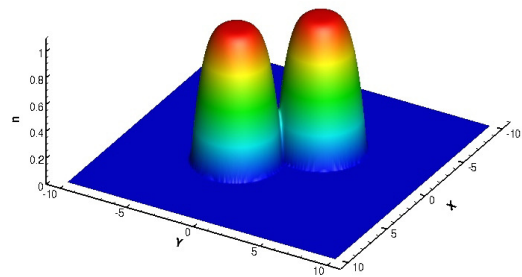}\label{fig:porous_5_t_01}}
\subfigure[t = 0.5]{\includegraphics[width=0.45\textwidth]{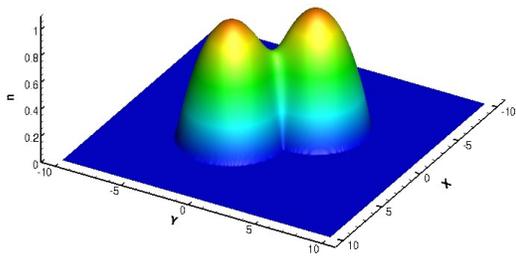}\label{fig:porous_5_t_05}}
\subfigure[t = 2]{\includegraphics[width=0.45\textwidth]{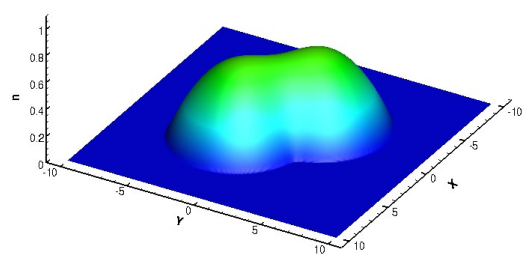}\label{fig:porous_5_t_20}}
\caption{Nonlinear porous medium problem with $h = 0.00733$.}
\label{fig:porous_time}
\end{figure}

\begin{figure}
\centering
\includegraphics[width=0.45\linewidth]{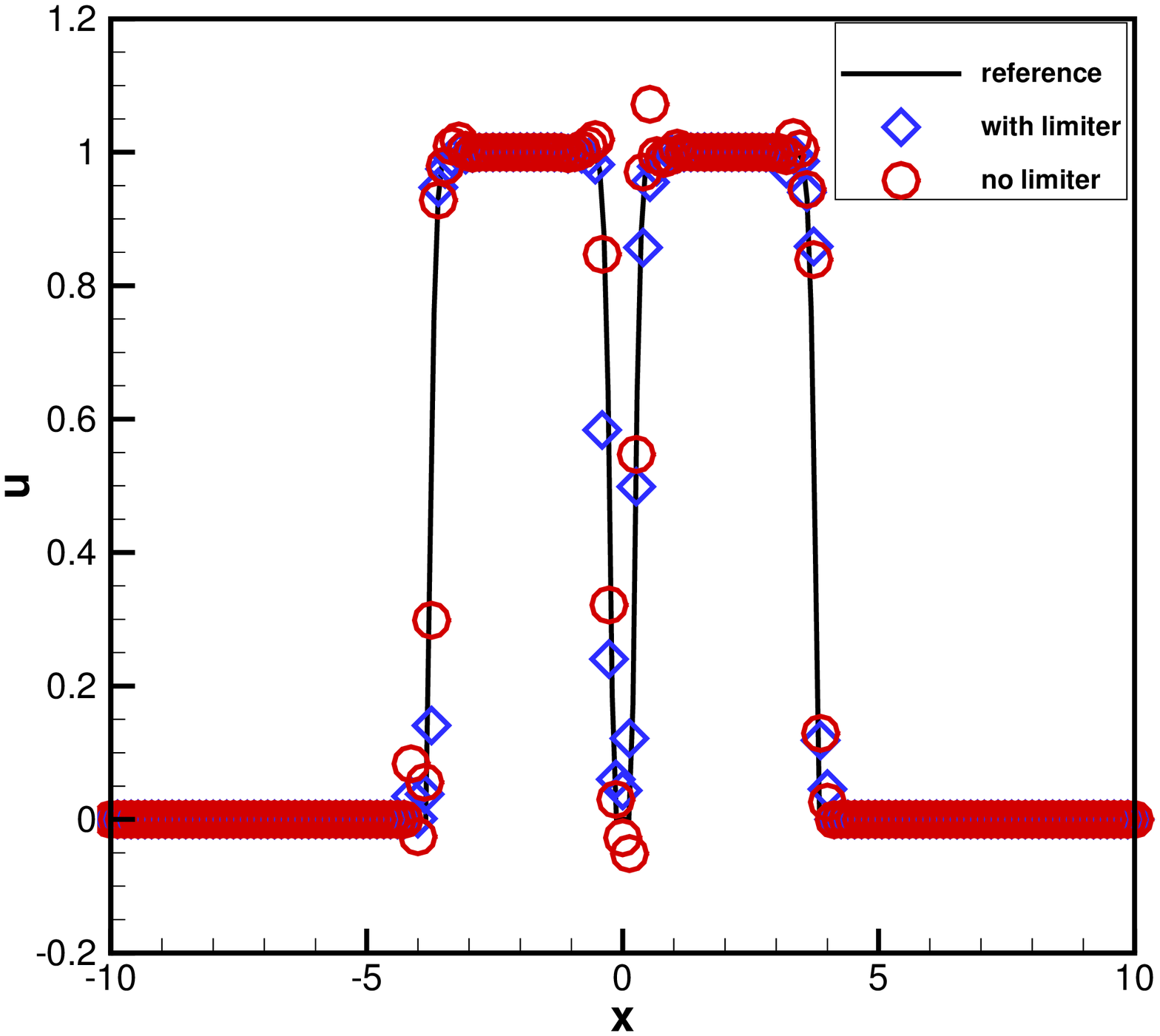}\includegraphics[width=0.45\linewidth]{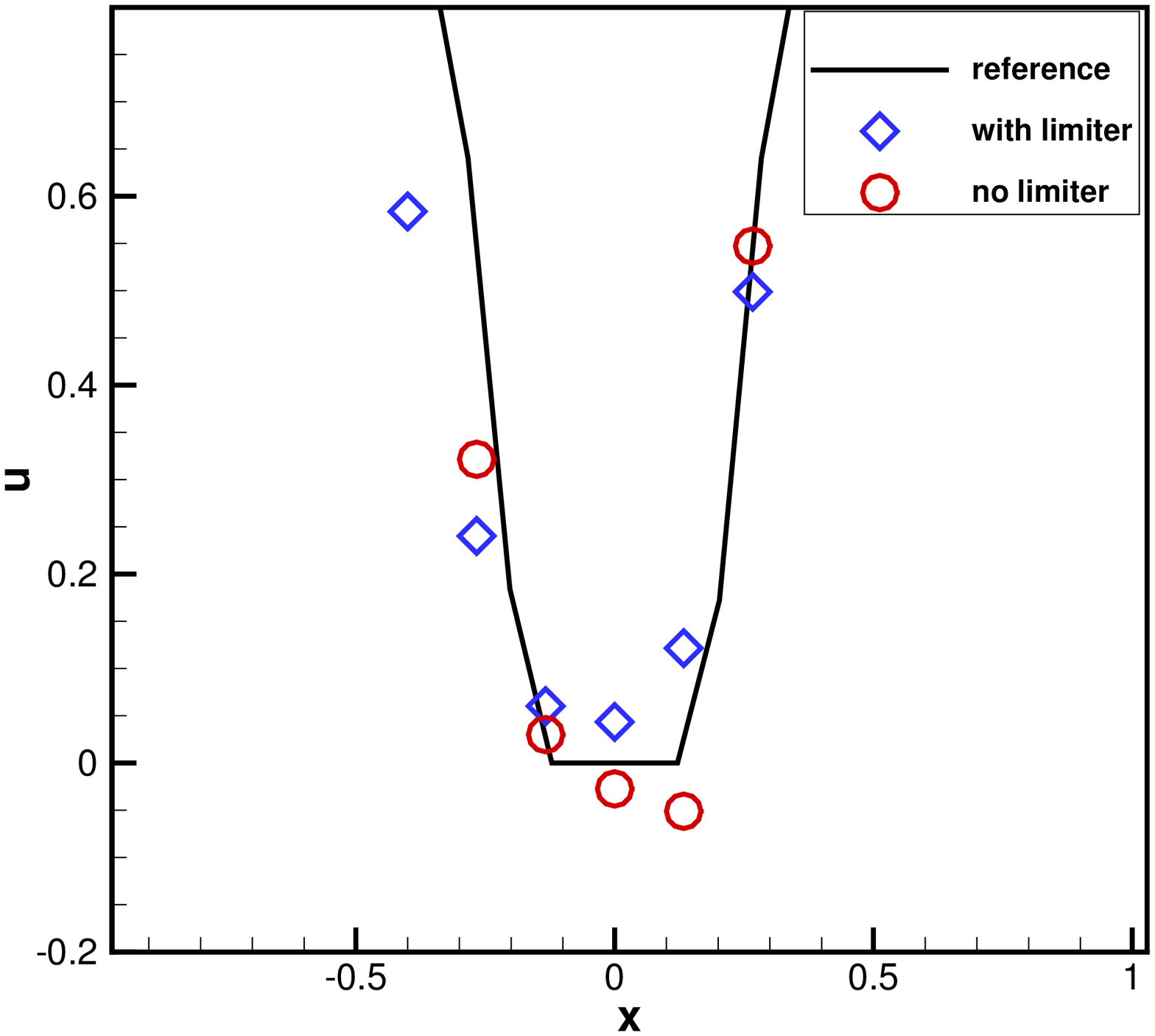}
\includegraphics[width=0.45\linewidth]{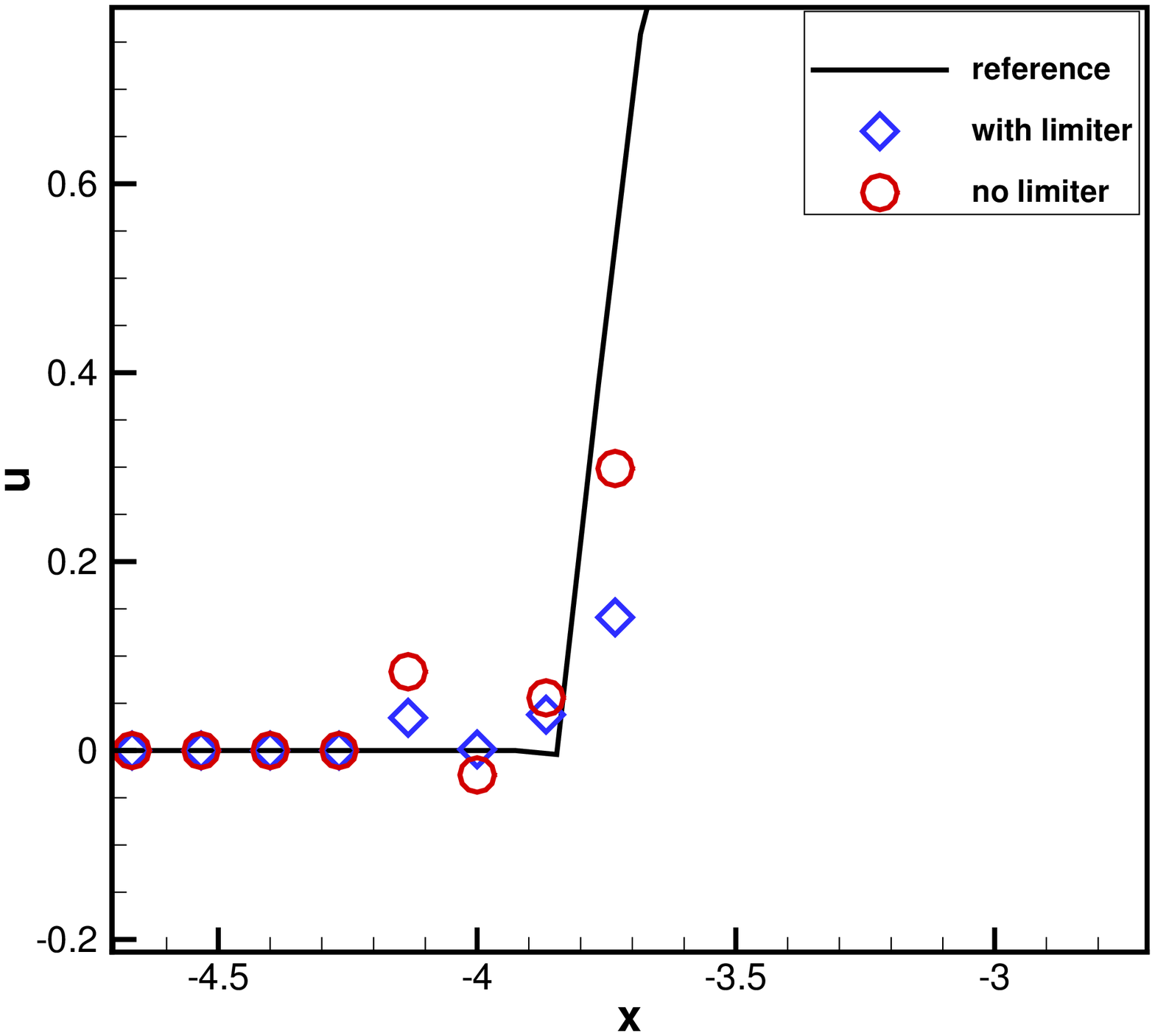}\includegraphics[width=0.45\linewidth]{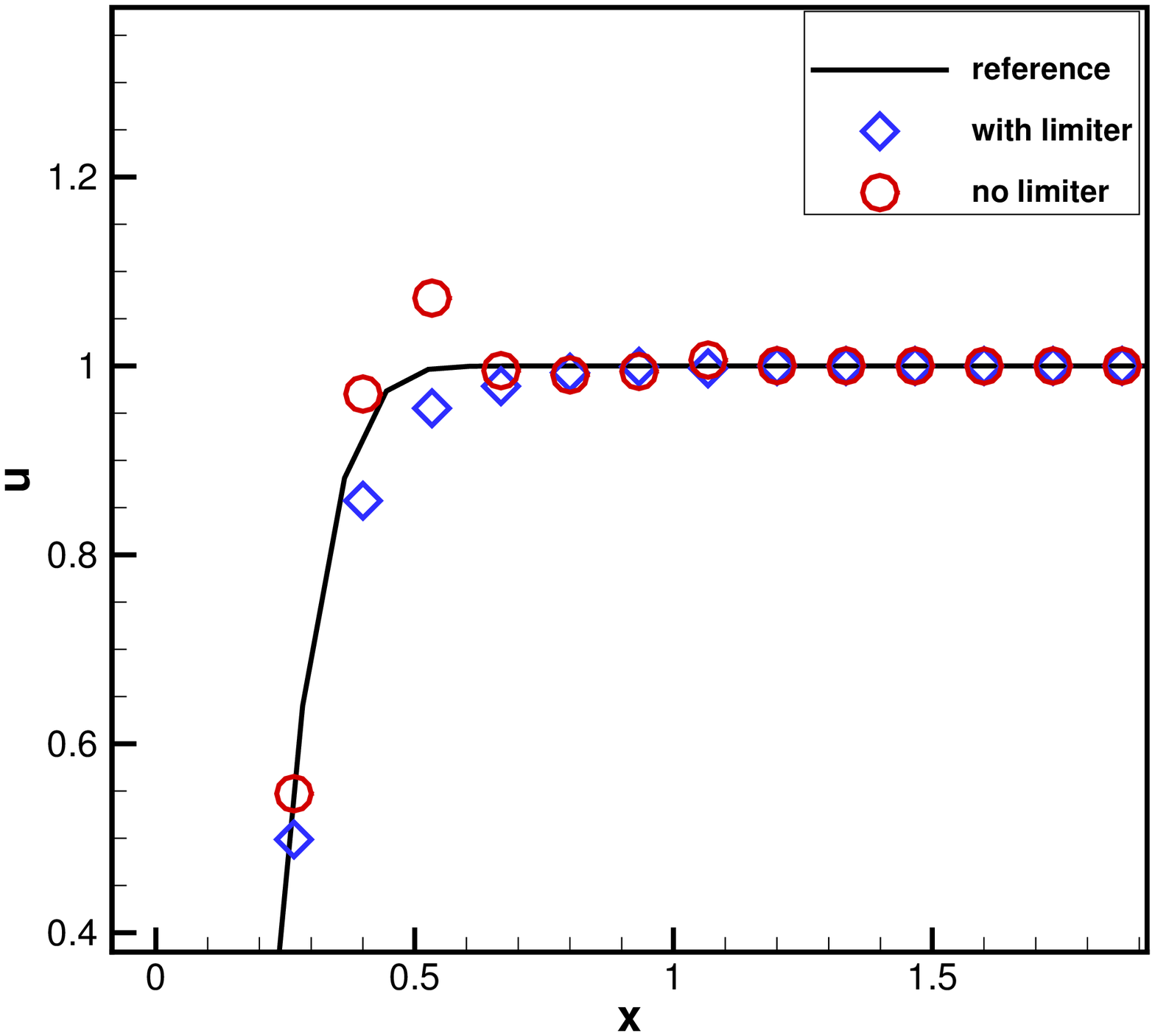}
\caption{The cut of the DDG interface correction solution along line $x + y = 0$ at $t = 0.005$.  Red circle symbol: no M-P-S limiter. Blue diamond symbol: M-P-S limiter applied.}
\label{fig:porous_compare}
\end{figure}

In this example we consider the nonlinear porous medium equation
\ben
u_t = \left(u^2\right)_{xx} +  \left(u^2\right)_{yy}, \quad (x,y) \in [-10,10]\times[-10,10],
\een
with initial condition
\ben
u(x,y,0) = \left\{
                \begin{array}{ll}
                  1,& (x-2)^2+(y+2)^2 < 6,\\
                  1,& (x+2)^2+(y-2)^2 < 6,\\
                  0,& otherwise,
                \end{array}
              \right.
\een
and zero boundary condition. Piecewise quadratic polynomial solutions implemented on unstructured mesh (Figure  \ref{fig:grid_3}) with size $h = 0.00733$ are shown in Figure \ref{fig:porous_time}. Notice that the minimum of the solution is zero. Numerical approximation without M-P-S limiter may become negative which may lead the problem ill-posed and cause the computations blow up. Implementations on a coarser mesh with $h = 0.0587$ are carried out, see Figure \ref{fig:porous_compare}. With M-P-S limiter applied, DDG interface correction solutions are maintained strictly inside the bound $[0,1]$.

\begin{figure}[htbp]
\centering
\includegraphics[width=0.65\linewidth]{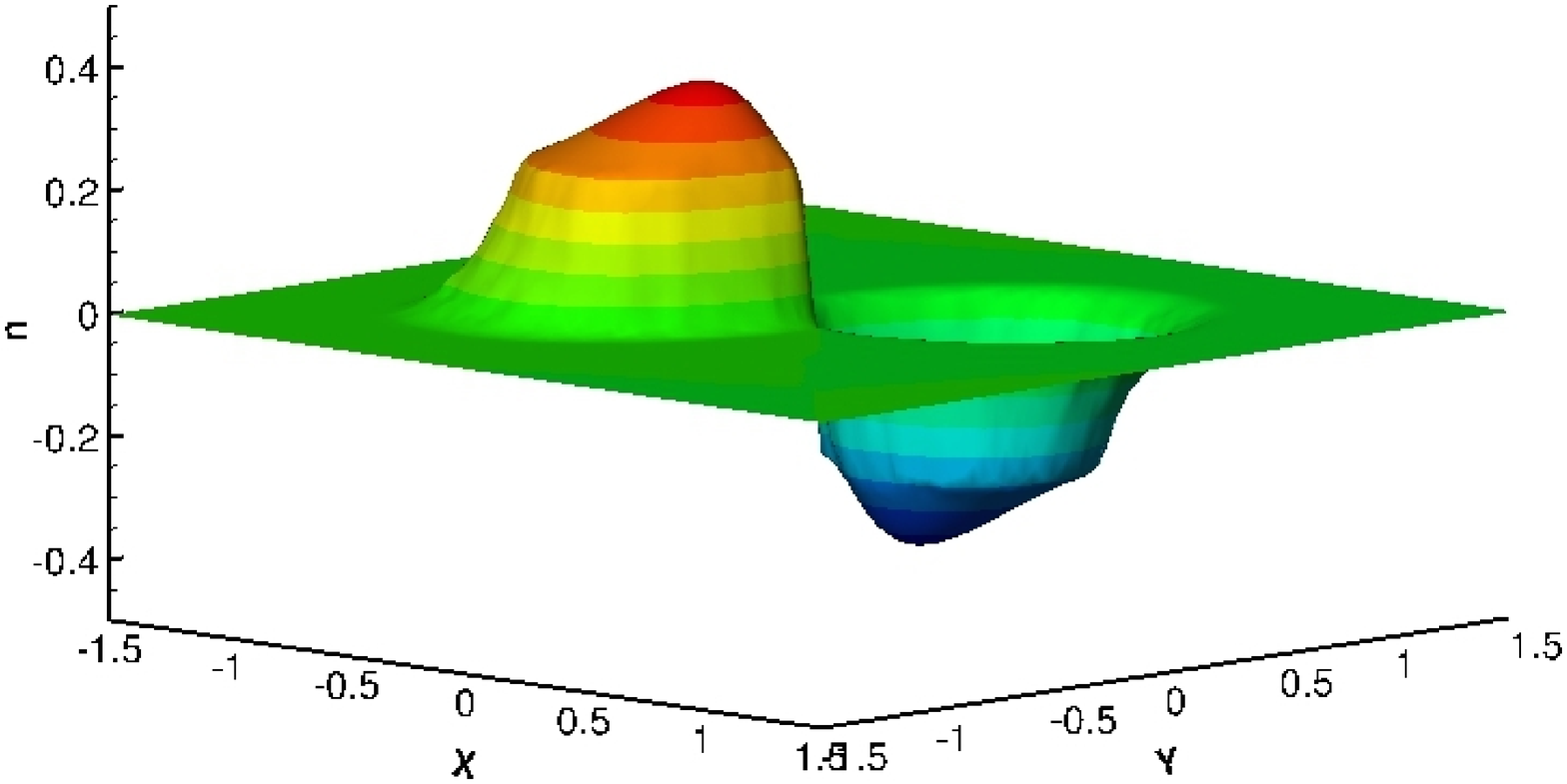}
\caption{Strongly degenerate parabolic problem solution at $t = 0.5$. }
\label{fig:strongly_degenerate_parabolic} 
\end{figure}

\begin{example}\label{e.g.:sdp} {Strongly degenerate parabolic problem}
\end{example}
We consider the following strongly degenerate parabolic problem with DDG interface correction method (\ref{scheme:conv-diff}),
\ben
\label{eq:sdp}
u_t + (u^2)_x + (u^2)_y = \epsilon(\nu(u)\nabla u_x)_x + \epsilon(\nu(u)\nabla u_y)_y, \quad (x,y) \in [-1.5,1.5]\times[-1.5,1.5].
\een
Initial condition is given with,
\ben
u(x,y,0) = \left\{
                \begin{array}{ll}
                   1,& (x+0.5)^2+(y+0.5)^2 < 0.16,\\
                  -1,& (x-0.5)^2+(y-0.5)^2 < 0.16,\\
                   0,& otherwise,
                \end{array}
              \right.
\een
and zero boundary condition is applied. We have $\epsilon = 0.1$ and,
\ben
\nu(u) = \left\{
                \begin{array}{ll}
                   0,& |u| \le 0.25,\\
                   1,& |u| > 0.25.
                \end{array}
              \right.
\een
Quadratic polynomial implementation with M-P-S limiter is carried out.  Here we apply the simple slope limiter \cite{Cockburn-book} to compress the oscillations caused from the nonlinear convection term. For the slope limiter, we take $\gamma = 1.5$ and $M = 5$. Implementation on mesh (Figure \ref{fig:grid_3}) with $h = 0.0147$ is shown in Figure \ref{fig:strongly_degenerate_parabolic}. The result agrees well those in literature, see \cite{Liu-Shu-Zhang2011, yan2013new}.

\begin{example}\label{e.g.:iNS} {Incompressible Navier-Stokes equation in vorticity stream-function formulation}
\end{example}

In this example, we consider to solve two-dimensional incompressible Navier-Stokes equation,
\be\label{eq:iNS}
\left\{
\begin{array}{ll}
&w_t + (uw)_x + (vw)_y = \frac{1}{Re}\nabla w,\\
&\Delta \phi = w, \qquad \langle u, v\rangle = \langle-\phi_y, \phi_x\rangle, \\
&\langle u, v\rangle \cdot {\bf n} = given, \qquad (x,y)\in \partial \Omega,
\end{array}
 \right.
\ee
written out in the vorticity stream-function format. We focus on the incompressible flow with high Reynolds numbers ($Re\gg 1$), thus explicit treatment on both convection term and diffusion term is efficient.

The initial vorticity is given with $w(x,y,0)=w_0(x,y)$ and periodic boundary condition is applied. We have $\phi(x,y)$ denoting the stream function and the velocity field is denoted as $\langle u, v\rangle$. We adopt the method of \cite{liu2000high} by Liu and Shu to solve (\ref{eq:iNS}). Thus at each time step, we first apply $P^2$ continuous finite element method as the Poisson solver to obtain stream function $\phi$, then have the velocity field and plug them into the vorticity equation and discretize in space with DDG interface correction method (\ref{scheme:conv-diff}), finally update the vorticity DG solution to the next time level. High order SSP Runge-Kutta explicit scheme \cite{Shu-Osher-TVD-1988} is applied to march forward solution in time. As remarked in \cite{liu2000high} that there is a natural match between the vorticity DG solution and the stream function. The normal component of velocity field $\langle u, v\rangle\cdot {\bf n}$ is continuous across all triangle edges, thus DG implementation on the convection part of vorticity equation is straight forward.

We carry out two tests in this example. First one is for accuracy check with exact solution maximum and minimum available and being applied with M-P-S limiter at each time step. The second one is a vortex patch problem.

\vspace{.15in}
\noindent
{\bf \large {Accuracy Test}}
\vspace{.15in}

We solve (\ref{eq:iNS}) with $Re = 100$.  Initial condition is $w_0(x,y) = -2\sin(x)\sin(y)$ with $\Omega = [0,2\pi]\times[0,2\pi]$ and periodic boundary condition is applied. Exact solution is available with
$$w(x,y,t) = -2\sin(x)\sin(y)\exp{(-2t/Re)}.$$
Quadratic $P^2$ implementations are carried out on mesh Figure \ref{fig:grid_3} and on mesh Figure \ref{fig:grid_5} with obtuse triangles. Errors and orders are listed in Table \ref{tab:iNS_accuracy_3} and Table \ref{tab:iNS_accuracy_5}. Again we have $we_{max}$ and $we_{min}$ represent exact solution maximum and minimum values. We observe that the M-P-S limiter removes all overshoots and undershoots and still maintains the third order accuracy.

\begin{table}[htbp]

   \begin{minipage}{1.\textwidth} 
    \resizebox{\textwidth}{!}{
\begin{tabular}{|c|cccc|cc|cccc|cc|}
\hline
$h/2\pi$ & \multicolumn{6}{c|}{without limiter} & \multicolumn{6}{c|}{with limiter} \\
\cline{2-13}
& $L^2$ error & Order & $L^\infty$ error & Order & $w_{\min}-we_{\min}$ & $w_{\max}-we_{\max}$ & $L^2$ error & Order & $L^\infty$ error & Order & $w_{\min}-we_{\min}$ & $w_{\max}-we_{\max}$\\
 \hline
0.117	&	1.90e-03	&		&	1.85e-02	&		&	-8.17e-03	&	7.72e-03	&	1.98e-03	&		 &	1.88e-02	 &		 &	0	&	0	\\ \hline 
0.0587	&	2.68e-04	&	2.82	&	2.07e-03	&	3.16	&	-3.64e-04	&	6.41e-04	&	 2.68e-04	&	2.88	&	 2.07e-03	&	3.18	&	0	&	0	\\ \hline 
0.0293	&	3.62e-05	&	2.89	&	2.92e-04	&	2.82	&	-6.40e-05	&	-4.50e-06	&	 3.62e-05	&	2.89	&	 2.92e-04	&	2.82	&	0	&	-4.52e-06	\\ \hline 
0.0147	&	4.54e-06	&	3.00	&	2.56e-05	&	3.51	&	-1.29e-06	&	3.62e-06	&	 4.54e-06	&	3.00	&	 2.56e-05	&	3.51	&	0	&	0	\\ \hline 
0.00733	&	5.84e-07	&	2.96	&	4.05e-06	&	2.66	&	-1.53e-07	&	1.95e-07	&	 5.84e-07	&	2.96	&	 4.05e-06	&	2.66	&	0	&	0	\\ \hline 
\end{tabular}}
\caption{Accuracy check on mesh Figure \ref{fig:grid_3}, final time $t=0.1$.}
\label{tab:iNS_accuracy_3}
   \end{minipage}
\vspace{0.5cm}

   \begin{minipage}{1.\textwidth} 
    \resizebox{\textwidth}{!}{
\begin{tabular}{|c|cccc|cc|cccc|cc|}
\hline
$h/2\pi$ & \multicolumn{6}{c|}{without limiter} & \multicolumn{6}{c|}{with limiter} \\
\cline{2-13}
& $L^2$ error & Order & $L^\infty$ error & Order & $w_{\min}-we_{\min}$ & $w_{\max}-we_{\max}$ & $L^2$ error & Order & $L^\infty$ error & Order & $w_{\min}-we_{\min}$ & $w_{\max}-we_{\max}$\\
 \hline
0.148	&	1.31e-03	&		&	2.74e-02	&		&	-8.40e-04	&	-7.18e-05	&	1.31e-03	&		 &	2.74e-02	 &		 &	0	&	-7.18e-05	\\ \hline 
0.0741	&	1.72e-04	&	2.93	&	3.35e-03	&	3.03	&	3.56e-05	&	-6.40e-05	&	 1.72e-04	&	2.93	&	 3.35e-03	&	3.03	&	4.15e-05	&	-6.40e-05	\\ \hline 
0.0371	&	2.37e-05	&	2.86	&	4.47e-04	&	2.90	&	3.36e-06	&	1.71e-06	&	 2.37e-05	&	2.86	&	 4.47e-04	&	2.90	&	3.36e-06	&	0	\\ \hline 
0.0185	&	2.99e-06	&	2.99	&	5.33e-05	&	3.07	&	1.79e-07	&	1.43e-07	&	 2.99e-06	&	2.99	&	 5.33e-05	&	3.07	&	1.79e-07	&	0	\\ \hline 
\end{tabular}}
\caption{Accuracy check on mesh Figure \ref{fig:grid_5} with obtuse triangles, final time $t=0.1$.}
\label{tab:iNS_accuracy_5}
   \end{minipage}
\end{table}

\vspace{.15in}
\noindent
{\bf \large {Vortex patch problem}}
\vspace{.15in}

\begin{figure}[htbp]
\centering
\includegraphics[width=0.45\linewidth]{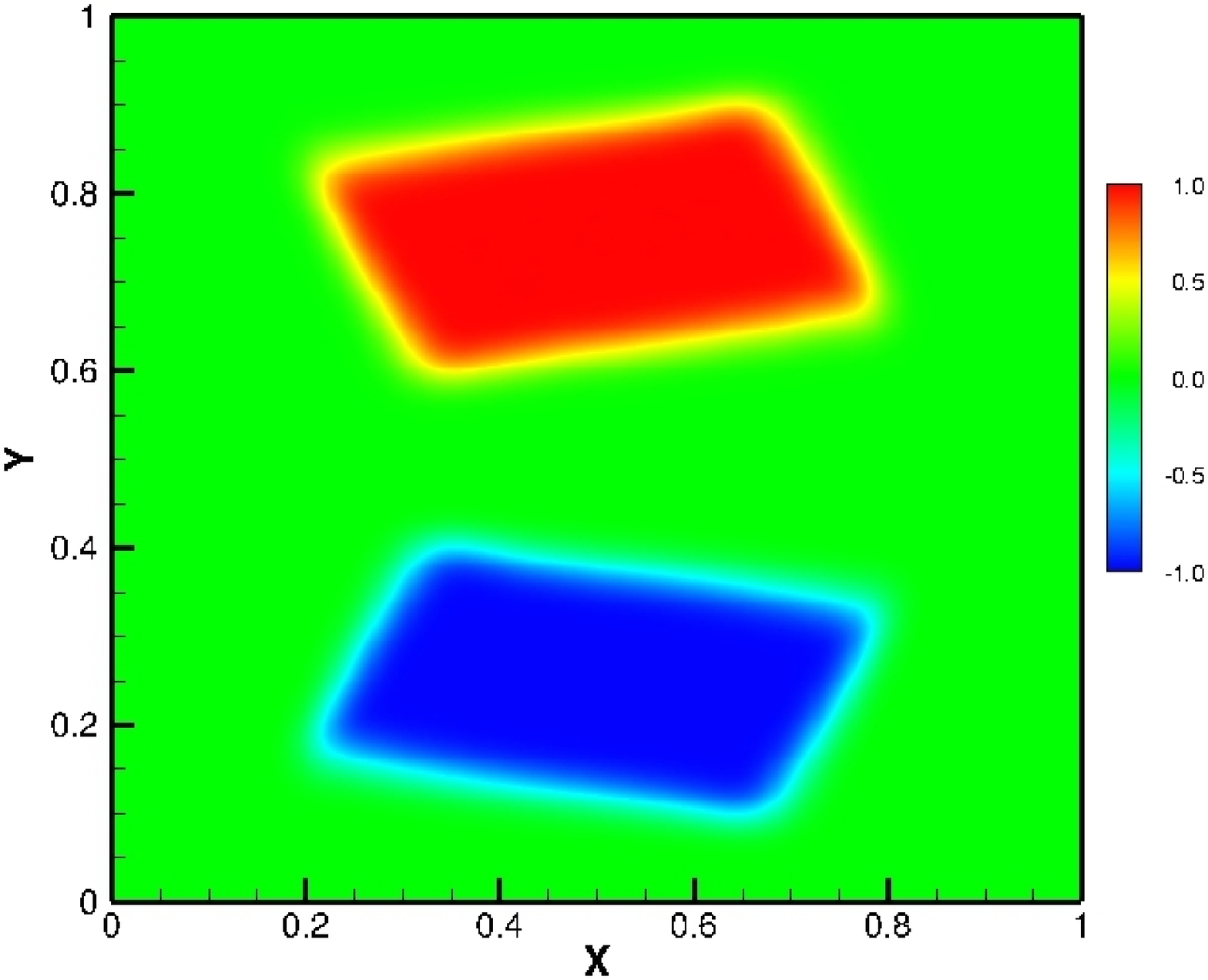}\includegraphics[width=0.45\linewidth]{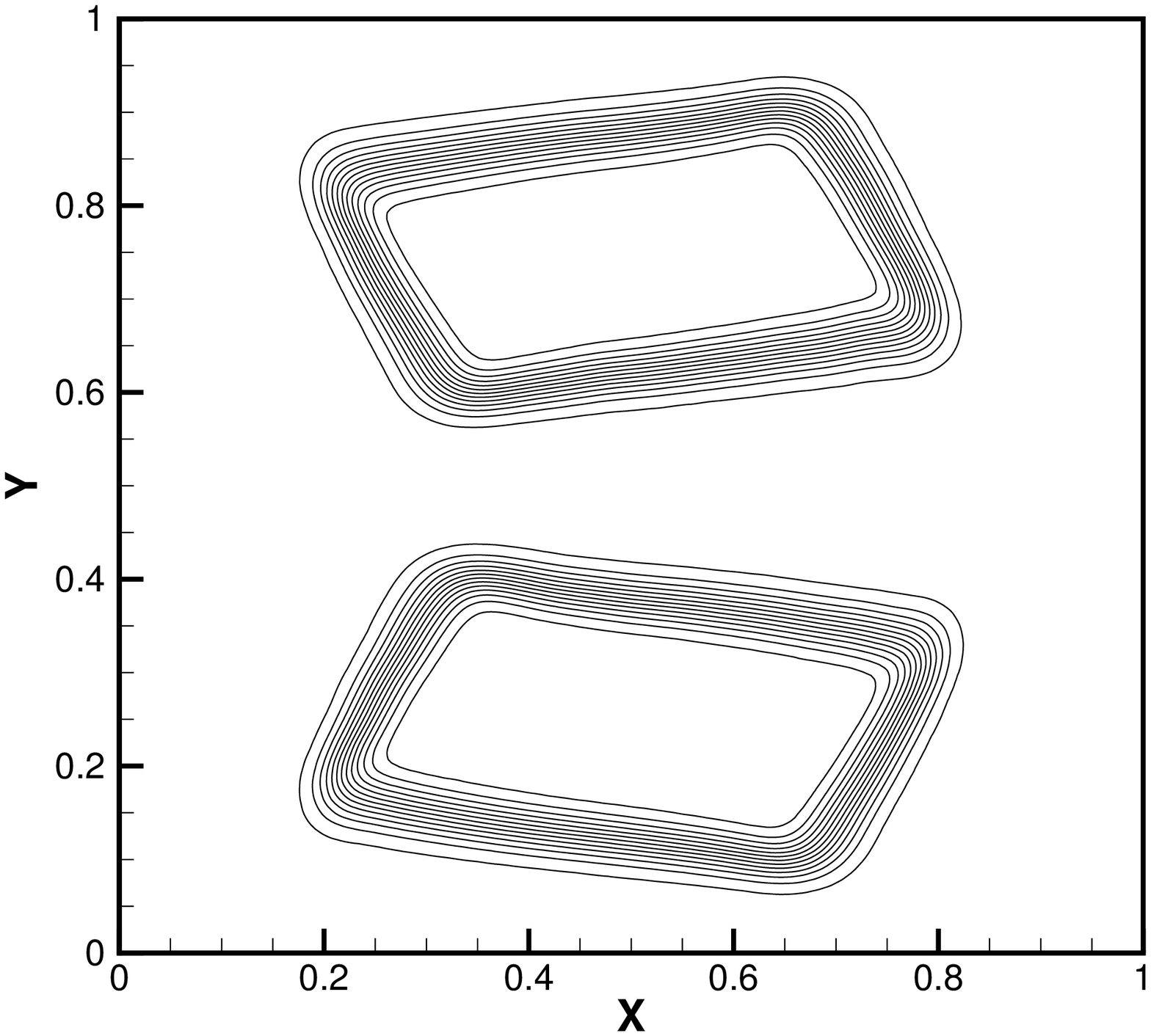}
\includegraphics[width=0.45\linewidth]{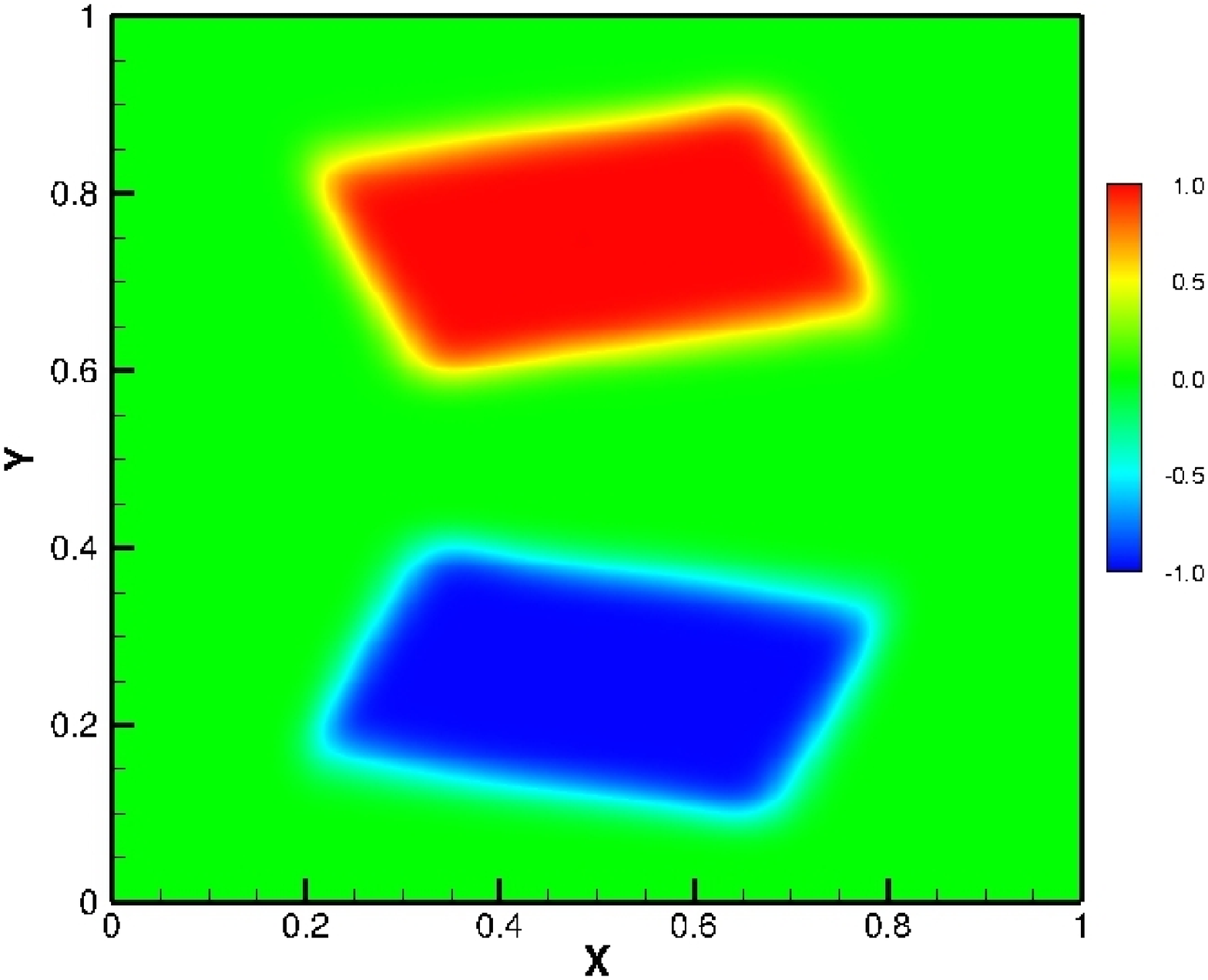}\includegraphics[width=0.45\linewidth]{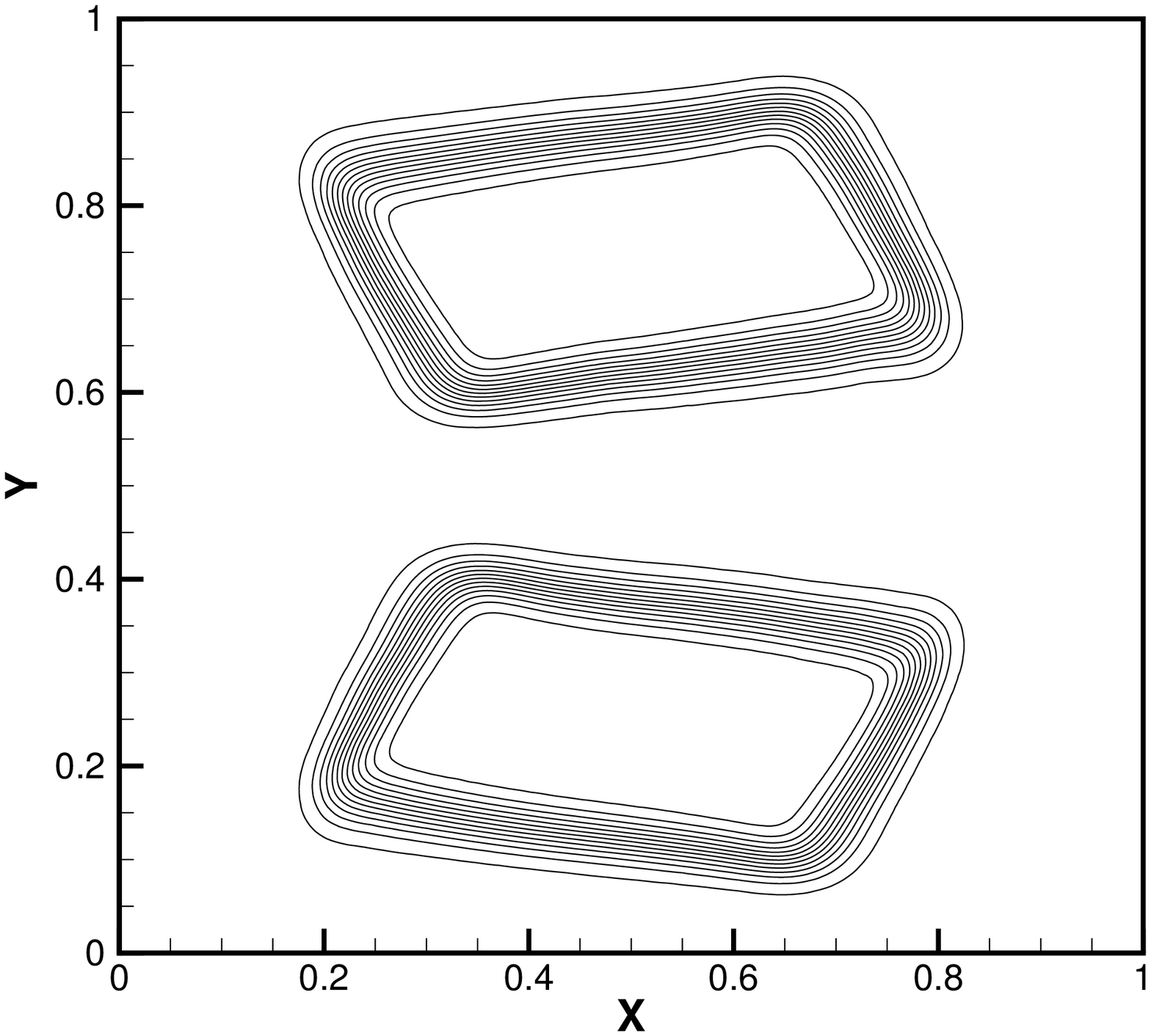}
\caption{Contours of the solutions, $Re = 100$, at $t = 1$ with mesh size $h = 0.0293\times2\pi$. Top: no M-P-S limiter ($w_{\min} =  -1.0011$, $w_{\max} = 1.0008$); Bottom: add M-P-S limiter ($w_{\min} =  -1$, $w_{\max} = 1$). 30 equally spaced contour lines are plotted.} 
\label{fig:iNS_vp_1_Re100_mesh3}
\end{figure}

\begin{figure}[htbp]
\centering
\includegraphics[width=0.45\linewidth]{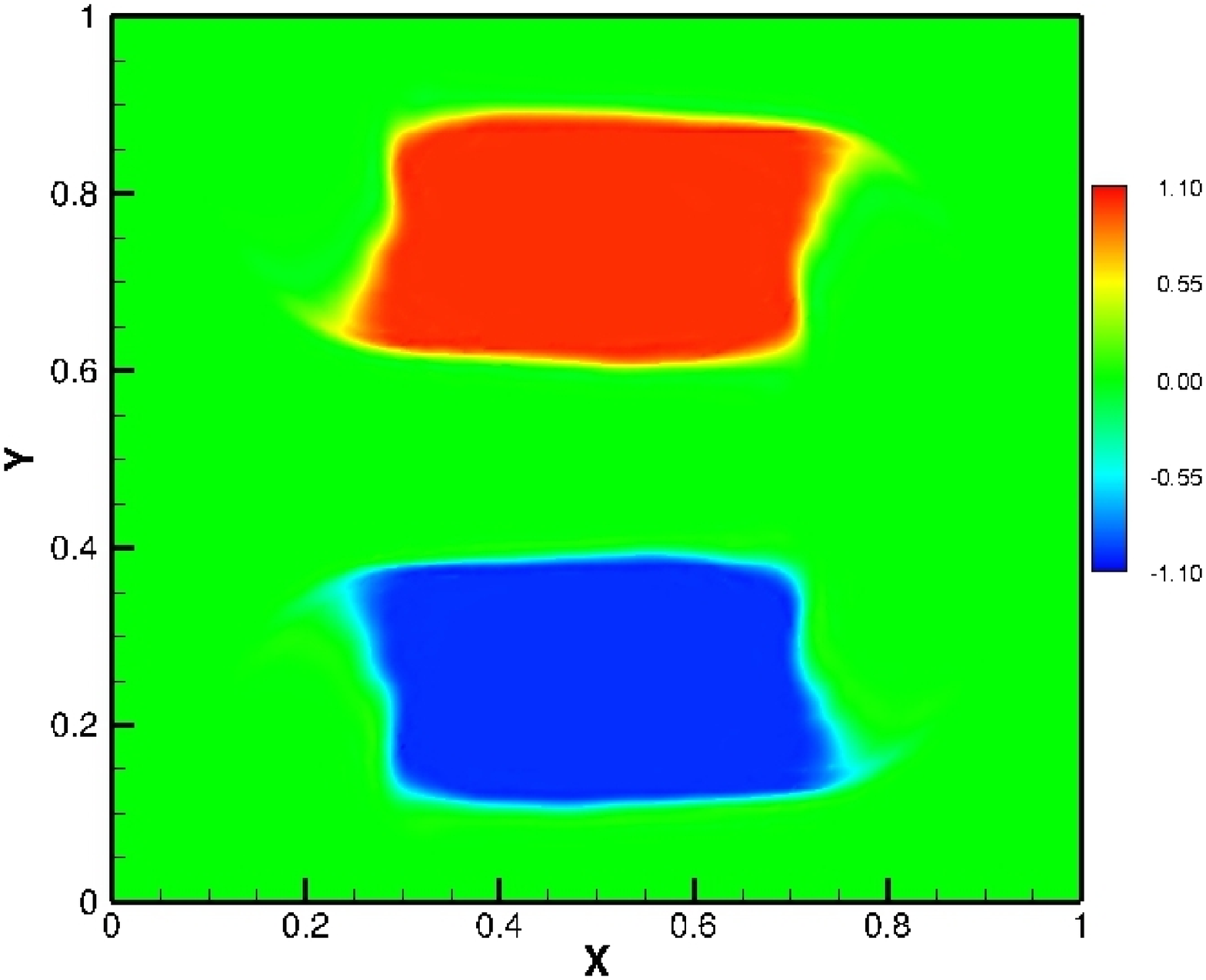}\includegraphics[width=0.45\linewidth]{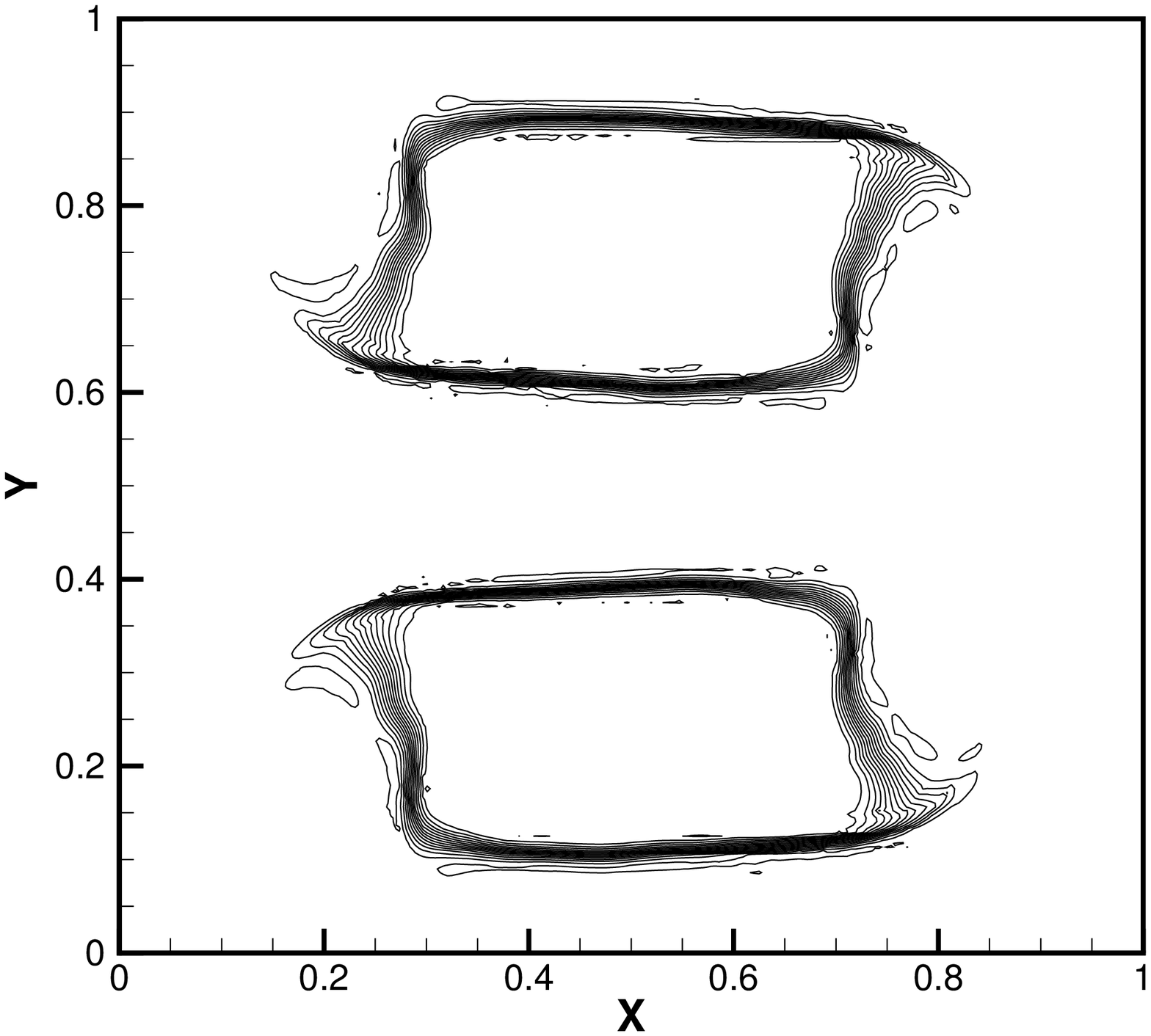}
\includegraphics[width=0.45\linewidth]{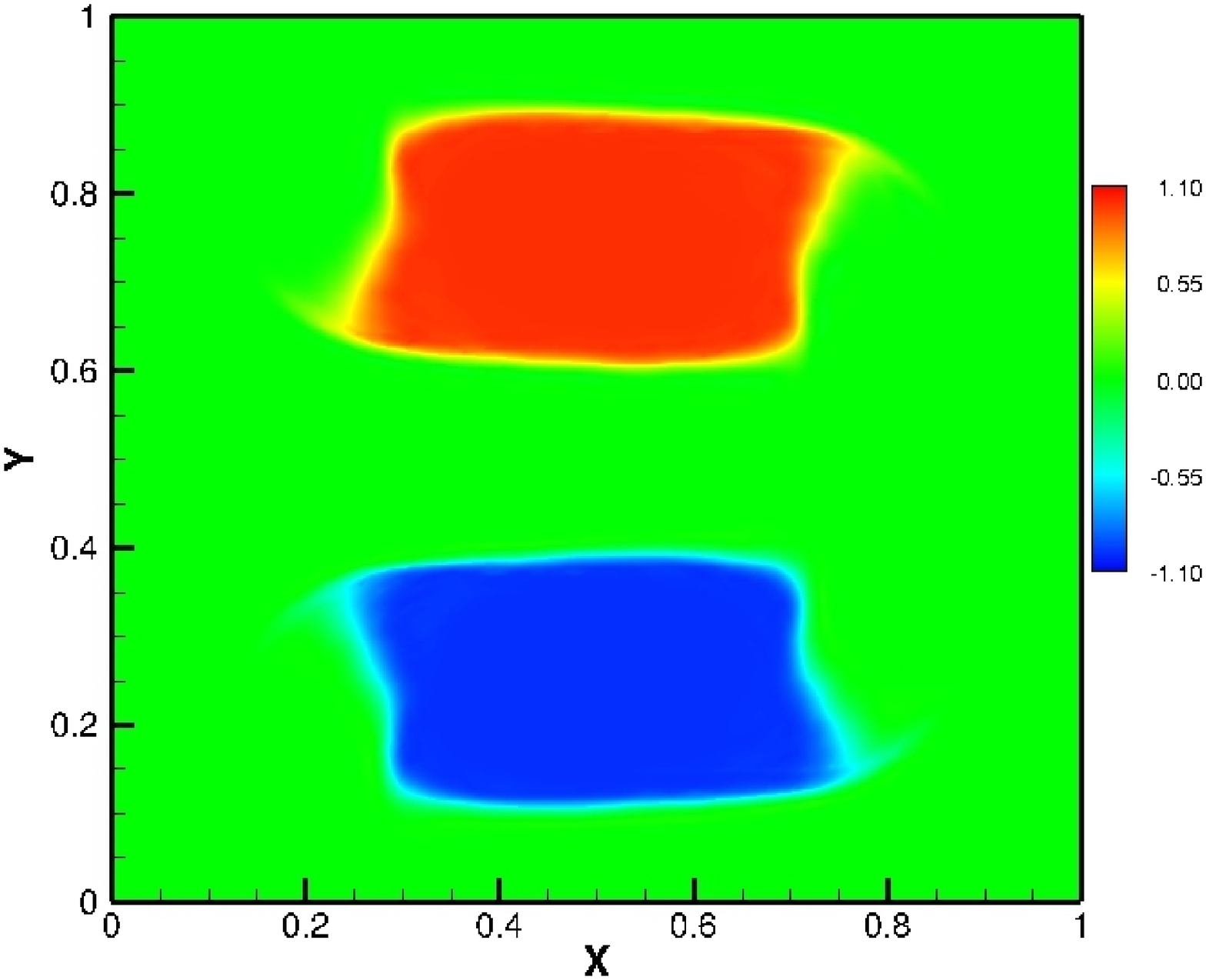}\includegraphics[width=0.45\linewidth]{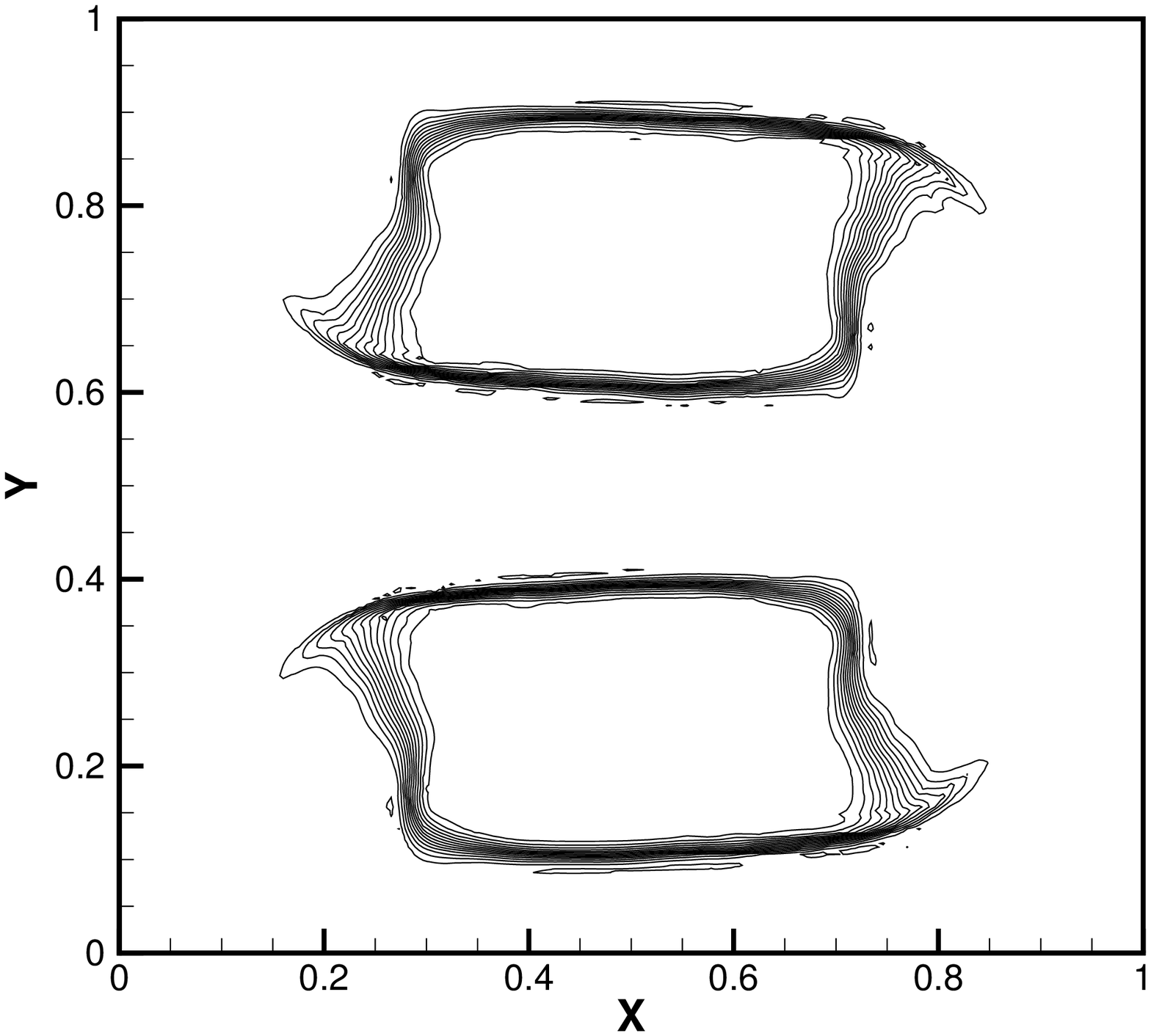}
\caption{Contours of the solutions, $Re = 10000$ at $t = 5$ and mesh size $h = 0.0293\times2\pi$. Top: no M-P-S limiter ($w_{\min} =  -1.1211$, $w_{\max} = 1.1508$); Bottom: add M-P-S limiter ($w_{\min} =  -1$, $w_{\max} = 1$). 30 equally spaced contour lines are plotted.} 
\label{fig:iNS_vp_1_Re10000_mesh3}
\end{figure}

Now we consider problem (\ref{eq:iNS})  with initial condition,
\be
w_0(x,y) =  \left\{
                \begin{array}{ll}
                  -1,& (x,y) \in [\frac{\pi}{2},\frac{3\pi}{2}] \times [\frac{\pi}{4},\frac{3\pi}{4}],\\
                  1,& (x,y) \in [\frac{\pi}{2},\frac{3\pi}{2}] \times [\frac{5\pi}{4},\frac{7\pi}{4}],\\
                  0,& otherwise,
                \end{array}
              \right.
\ee
and periodic boundary condition. The Reynolds number is chosen to be $Re = 100$ or $Re = 10000$. We compare the maximum and minimum of the numerical solutions with and without M-P-S limiter applied, see Table \ref{tab:iNS_vp_0.1_Re100_mesh1} and Table \ref{tab:iNS_vp_0.1_Re10000_mesh1}. Mesh Figure \ref{fig:grid_3} is used and quadratic polynomials is applied. We also plot the solutions for the case $Re = 100$ at $t = 1$, shown in Figure \ref{fig:iNS_vp_1_Re100_mesh3}, and the case $Re = 10000$ at $t = 5$, shown in Figure \ref{fig:iNS_vp_1_Re10000_mesh3}. It is clear that the overshoots and undershoots are removed with the M-P-S limiter applied.

\begin{table}
\centering
\begin{tabular}{|c|cc|cc|}
\hline
$Re = 100$  & \multicolumn{2}{c|}{without limiter} & \multicolumn{2}{c|}{with limiter} \\
\hline
 $h/2\pi$& $w_{\min}-we_{\min}$ & $w_{\max}-we_{\max}$ & $w_{\min}-we_{\min}$ & $w_{\max}-we_{\max}$\\
 \hline
0.117	&	-6.96e-01	&	6.87e-01	&	0	&	0	\\ \hline 
0.0587	&	-2.72e-01	&	2.66e-01	&	0	&	0	\\ \hline 
0.0293	&	-8.02e-02	&	4.22e-02	&	0	&	0	\\ \hline 
0.0147	&	-3.17e-03	&	3.00e-03	&	0	&	0	\\ \hline 
0.00733	&	-6.41e-04	&	7.41e-04	&	0	&	0	\\ \hline 
\end{tabular}
\caption{Maximum and minimum of the solutions, $Re = 100$ at $t = 0.1$.}
\label{tab:iNS_vp_0.1_Re100_mesh1}
\vspace{0.3cm}
\begin{tabular}{|c|cc|cc|}
\hline
$Re = 10000$  & \multicolumn{2}{c|}{without limiter} & \multicolumn{2}{c|}{with limiter} \\
\hline
 $h/2\pi$& $w_{\min}-we_{\min}$ & $w_{\max}-we_{\max}$ & $w_{\min}-we_{\min}$ & $w_{\max}-we_{\max}$\\
 \hline
0.117	&	-8.62e-01	&	8.62e-01	&	0	&	0	\\ \hline 
0.0587	&	-6.14e-01	&	7.82e-01	&	0	&	0	\\ \hline 
0.0293	&	-5.14e-01	&	5.77e-01	&	0	&	0	\\ \hline 
0.0147	&	-4.35e-01	&	4.79e-01	&	0	&	0	\\ \hline 
0.00733	&	-3.52e-01	&	2.81e-01	&	0	&	0	\\ \hline  
\end{tabular}
\caption{Maximum and minimum of the solutions, $Re = 10000$ at $t = 0.1$.}
\label{tab:iNS_vp_0.1_Re10000_mesh1}

\end{table}


%

\section{Appendix}\label{sec:appendix}
\numberwithin{equation}{section}

As discussed in Theorem \ref{thm:linear}, we need to find a 13-point quadrature rule that is exact for $P^2$ polynomial and include all the selected points (total 12) as quadrature points. For example we consider the edge $AB$ shared by $K$ and $K_3$. As shown in Figure \ref{fig:selectedpoints}, five points ${\bf x}_{K,i}$ ($i = 1,\cdots,5$) from element $K$'s side are used in the calculation of (\ref{h_{AB}}) - (\ref{un_edge_AB}). In this section, our goal is to find such a quadrature rule with the selected points as quadrature points and having non-negative weights. Specifically we need to find the minimum weight in the quadrature rule, since it is used in (\ref{cond}) - (\ref{A_cfl}) to bound CFL condition and further identify suitable time step size for time evolution. For convenience, we introduce notation $|K|$ to represent the area of triangle element $K$. 

\subsection{Quadrature rules for cell average}\label{quadrature}
Again we use edge $AB$ to illustrate the way to find such quadrature rule with ${\bf x}_{K,i}$ ($i = 1,\cdots,5$) as quadrature points. The method we investigate is that we look for weights $\omega_i$ in the following format,

$$
\frac{1}{3}\overline{u}^n_K=\frac{1}{3|K|}\int_Ku^n_K(x,y) dxdy=\sum_{i=1}^5 w_i u^n_K({\bf x}_{K,i})+\sum_{j=1}^{l}w^*_j u^n_K({\bf x}_{K,j}^*),
$$
with other quadrature points ${\bf x}_{K,j}^*$ ($j=1,\cdots,l$) selected in the following to complete the quadrature rule.  Notice we will combine these points together with other points from edges $BC$ and $CA$ to find the location and weight for the 13-th point.

\begin{figure}[htbp]
\centering
\subfigure[Angles in $K$ and $K_3$]{\includegraphics[width=0.32\linewidth]{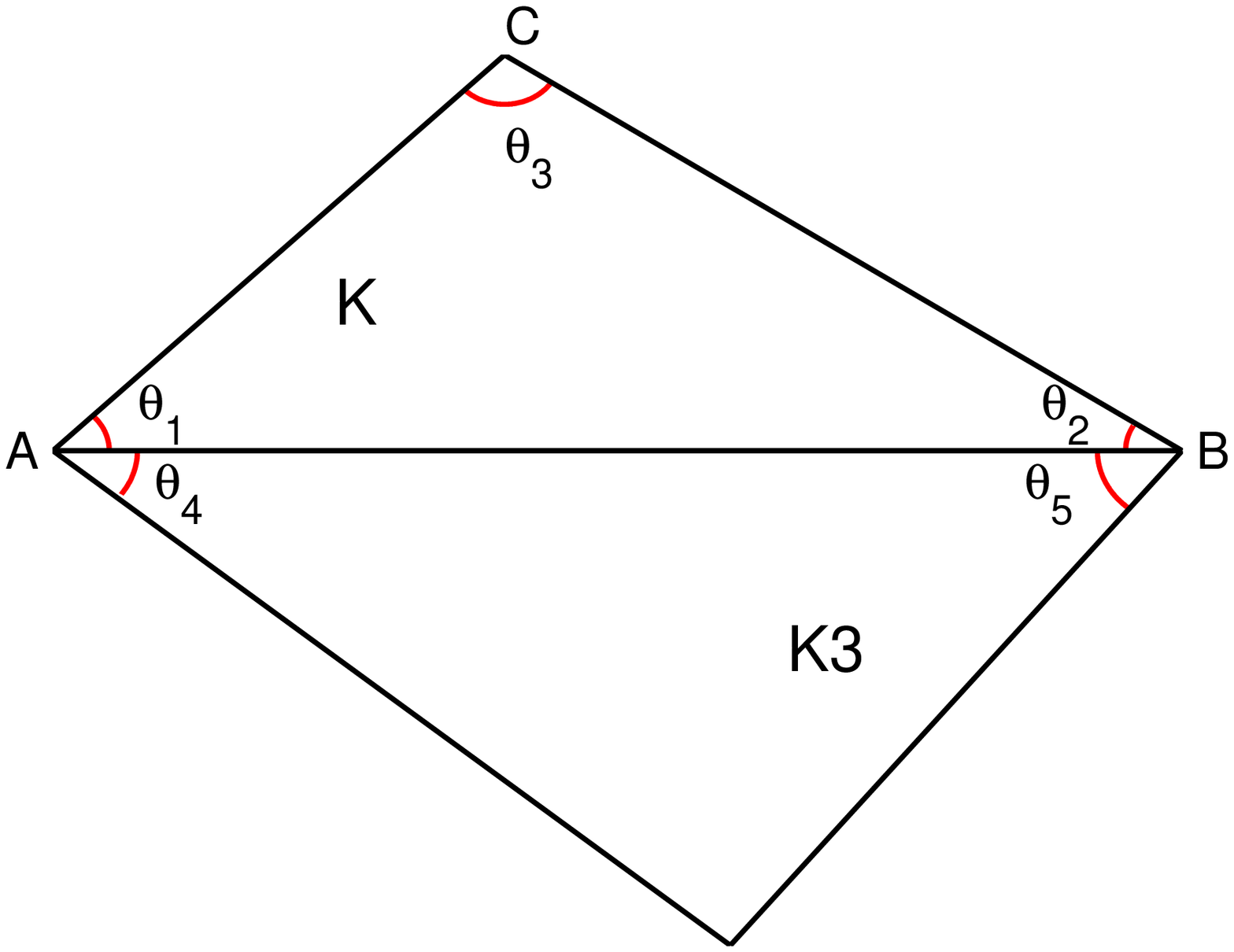}
\label{fig:neighbourtriangles}}
\subfigure[Case 1]{\includegraphics[width=0.32\textwidth]{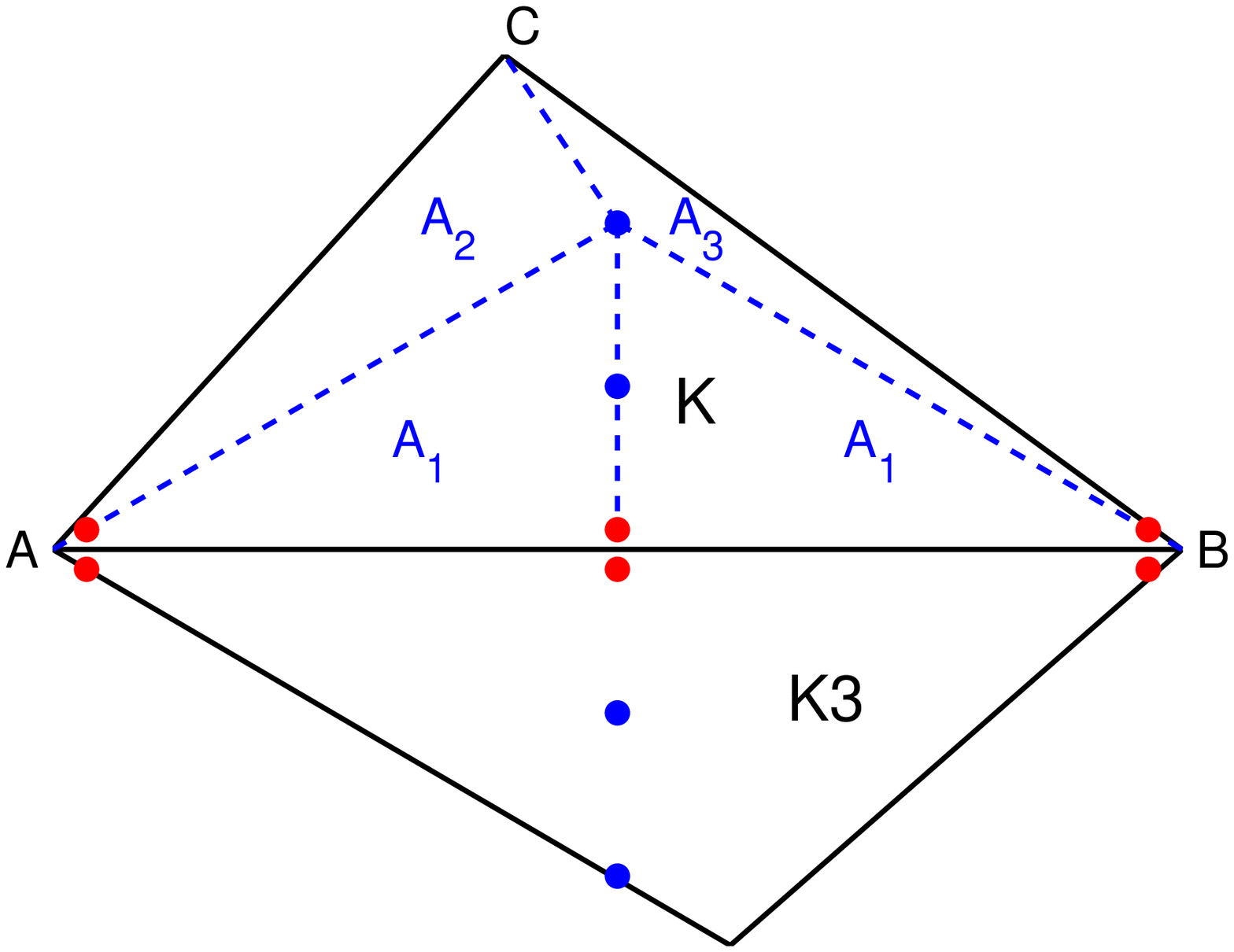}\label{fig:case1}}
\subfigure[Case 2]{\includegraphics[width=0.32\textwidth]{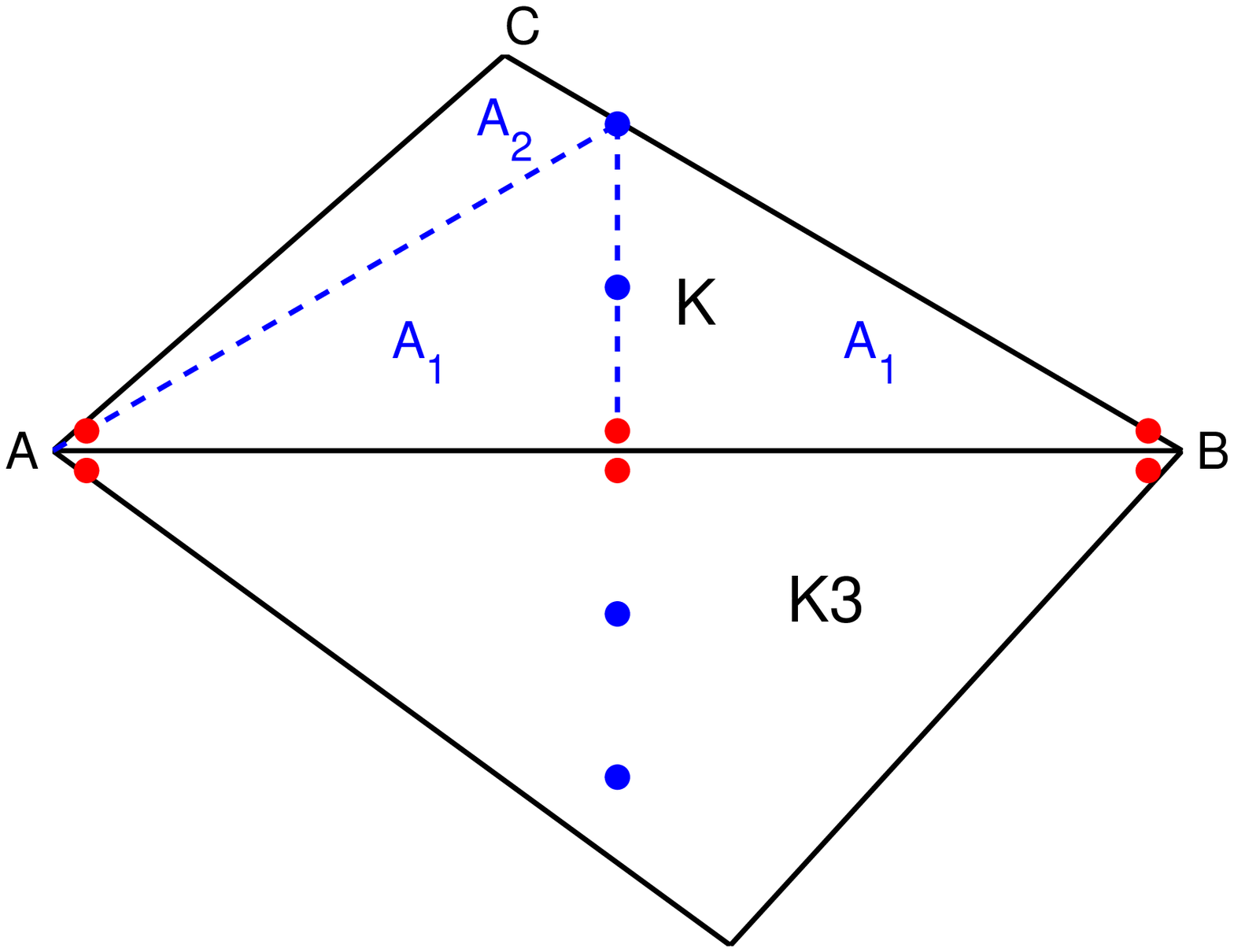}\label{fig:case2}}
\caption{Left: Element $K$  with its neighbor $K_3$; Right: Two cases of points selected along normal vector($A1$, $A2$, $A3$ denote small triangles area).}
\label{fig:triangle}
\end{figure}

The subsection \ref{quadGaussL-GaussR} below offers one way to find a quadrature rule on triangle element with vortices included as quadrature points (with weight $\omega_1$ in (\ref{eq:weight4vertex})). We also use a quadrature rule for triangle element only with edge centers with equal weights $1$. We divide triangle $K$ into small triangles with ${\bf x}_{K,i}$ ($i = 1,\cdots,5$) being vertices or edge centers and use each rule for a half $\frac{1}{6}\overline{u}^n_K$. In this way, all other vertices and edge centers are selected as quadrature points as well.

Let's focus the discussion on point ${\bf x}_{K,4}$ and ${\bf x}_{K,5}$, shown in Figure \ref{fig:selectedpoints} or the blue dots in Figure \ref{fig:triangle}.
First of all, the location of ${\bf x}_{K,4}$ and ${\bf x}_{K,5}$ are determined by elements $K$ and $K_3$. There are two cases: $(1)$ $\min(\theta_1,\theta_2)>\min(\theta_4,\theta_5)$, then ${\bf x}_{K,5}$ is inside $K$, see Figure \ref{fig:case1};  $(2)$ $\min(\theta_1,\theta_2)\le\min(\theta_4,\theta_5)$,  then ${\bf x}_{K,5}$ is on the edge of $K$, see Figure \ref{fig:case2}. The triangle geometrical information is shown and marked in Figure \ref{fig:triangle}.  We can estimate the quadrature weights of ${\bf x}_{K,i}$ ($i = 1,\cdots,5$) in each case, see Table \ref{Estimate on the quadrature weights}.

\begin{table}
\begin{center}
\begin{tabular}{|c||c|c|}
\hline
Point & Case 1 & Case 2 \\ \hline \hline
${\bf x}_{K,4}$
& $ \frac{1}{6}\cdot \left(\frac{2A_1}{|K|}\cdot \frac{1}{3}\right)
\geq \frac{1}{18}\tan(\hat{\theta})\cot(\check{\theta})$
& $ \frac{1}{6}\cdot \left(\frac{2A_1}{|K|}\cdot \frac{1}{3}\right)
\geq \frac{1}{18}\tan(\hat{\theta})\cot(\check{\theta})$
\\ \hline
${\bf x}_{K,1}$
& $\frac{1}{6}\cdot \left(\frac{A_1+A_2}{|K|}w_1\right)
\geq \frac{w_1}{12}\tan(\hat{\theta})\cot(\check{\theta})$
& $\frac{1}{6}\cdot \left(\frac{A_1+A_2}{|K|}w_1\right)
\geq \frac{w_1}{12}\tan(\hat{\theta})\cot(\check{\theta})$
\\ \hline
${\bf x}_{K,2}$
& $\frac{1}{6}\cdot \left(\frac{A_1+A_3}{|K|}w_1 \right)
\geq \frac{w_1}{12}\tan(\hat{\theta})\cot(\check{\theta})$
& $\frac{1}{6}\cdot \left( \frac{A_1}{|K|}w_1\right)
\geq \frac{w_1}{12}\tan(\hat{\theta})\cot(\check{\theta})$
\\ \hline
${\bf x}_{K,3}$
& $\frac{1}{6}\cdot \left(\frac{2A_1}{|K|}w_1 \right)
\geq \frac{w_1}{6}\tan(\hat{\theta})\cot(\check{\theta})$
& $\frac{1}{6}\cdot \left( \frac{2A_1}{|K|}w_1\right)
\geq \frac{w_1}{6}\tan(\hat{\theta})\cot(\check{\theta})$
\\ \hline
${\bf x}_{K,5}$
& $\frac{1}{6} \cdot \left( \frac{2A_1+A_2+A_3}{|K|}w_1 \right) = \frac{w_1}{6}$
& $\frac{1}{6} \cdot \left( \frac{2A_1+A_2}{|K|}w_1\right)  = \frac{w_1}{6}$\\ \hline
\end{tabular}
\end{center}
\caption{Estimate of the quadrature weights}
\label{Estimate on the quadrature weights}
\end{table}%
\vspace{.2in}

We apply the same procedure to edge $BC$ and $CA$ to include selected points as quadrature points. Collect all data from the three edges, we have the estimate on the weights as follows,
\be\label{weights}
\frac{\partial \overline{u}^n_K}{\partial u_{K,1}}, \frac{\partial \overline{u}^n_K}{\partial u_{K,2}}, \frac{\partial \overline{u}^n_K}{\partial u_{K,3}}
\geq \frac{w_1}{6}\tan(\hat{\theta})\cot(\check{\theta}), \qquad
\frac{\partial \overline{u}^n_K}{\partial u_{K,4}}
\geq  \frac{1}{18}\tan(\hat{\theta})\cot(\check{\theta}), \qquad
\frac{\partial \overline{u}^n_K}{\partial u_{K,5}}
= \frac{w_1}{6}.
\ee

Now we are ready to find a non-negative quadrature rule with only 13 points inside triangle K which include the 12 selected points from numerical flux integral on the edges. Let's reorder the quadrature points, and denote $u_{K,i}$ ($i = 1, \cdots, 12$)  as the 12 points selected from evaluating $\int_{\partial K} \widehat{u_{\bf n}} ds$ with weights $w_i$ ($i = 1, \cdots, 12$). We also reorder the rest points and denote them as $u^*_{K,j}$ ($j = 1, \cdots, m$, where $m$ is an integer), with weights $w^*_j$ ($j = 1, \cdots, m$). Let $w_{13} = \sum_{j=1}^m w^*_j$, we have $\sum_{j=1}^{13} w_j = \sum_{j=1}^{12} w_j + \sum_{j=1}^m w^*_j = 1.$
Moreover,
$\frac{1}{w_{13}}\sum_{j=1}^m w^*_j u^*_{K,j} = \sum_{j=1}^m \frac{w^*_j}{w_{13}}u^*_{K,j}$
 is a convex combination of the points $u^*_{K,j}$ ($j = 1, \cdots, m$). By the mean value theory, one can find a point $u_{K,13}$ inside the convex hall of the points $u^*_{K,j}, j = 1, \cdots, m$ such that $u_{K,13} = \sum_{j=1}^m \frac{w^*_j}{w_{13}}u^*_{K,j}$. Therefore, we have a quadrature rule with positive weights of total 13 points inside triangle $K$ which include the 12 selected points: $\left\{ \left(u_{K,j}, w_j\right) | j = 1, \cdots, 12 \right\}$.

\subsection{Quadrature rule for triangle element with vertices}\label{quadGaussL-GaussR}
In this section, we design one quadrature rule that is exact for quadratic polynomial $P^2(K)$ on any triangle element $K$. Especially we include the three vortices and three edge centers in the quadrature points set. This work is inspired by the quadrature rule designed in \cite{zhangxxTri}. Specifically we are interested in finding the weights before the vertices.

For convenience, we use the position vectors to denote the three vertices of $K$: ${\bf v}^1$, ${\bf v}^2$ and ${\bf v}^3$. Thus, the position vector ${\bf P}$ of any point $P$ inside triangle $K$ can be described by the barycentric coordinates $(\xi_1, \xi_2, \xi_3)$, i.e., ${\bf P}=\xi_1{\bf v}^1+\xi_2{\bf v}^2+\xi_3{\bf v}^3$.

We first consider the quadrature rule on the unit square with vertices coordinates as $\left(-\frac{1}{2},-\frac{1}{2}\right)$, $\left(-\frac{1}{2},\frac{1}{2}\right)$, $\left(\frac{1}{2},-\frac{1}{2}\right)$ and $\left(\frac{1}{2},\frac{1}{2}\right)$ in the $u$-$v$ plane, and then we use projections/transformations to map it to a quadrature rule on triangle element $K$.

\begin{figure}[!htbp]
\centering
\subfigure[$S^2$]{\includegraphics[width=0.24\textwidth]{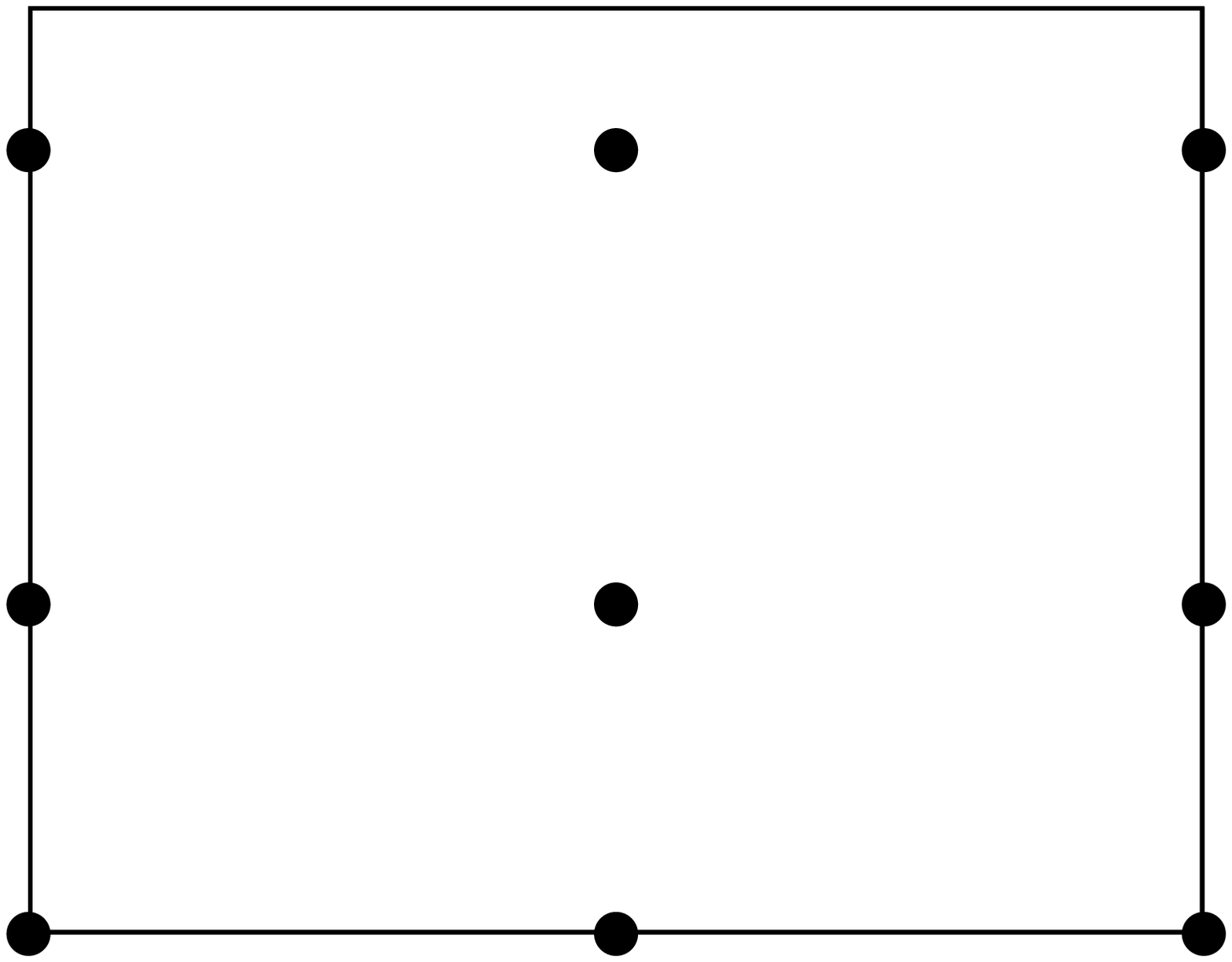}\label{fig:quadS2K_0}}
\subfigure[${\bf g}_1(S^2)$]{\includegraphics[width=0.24\textwidth]{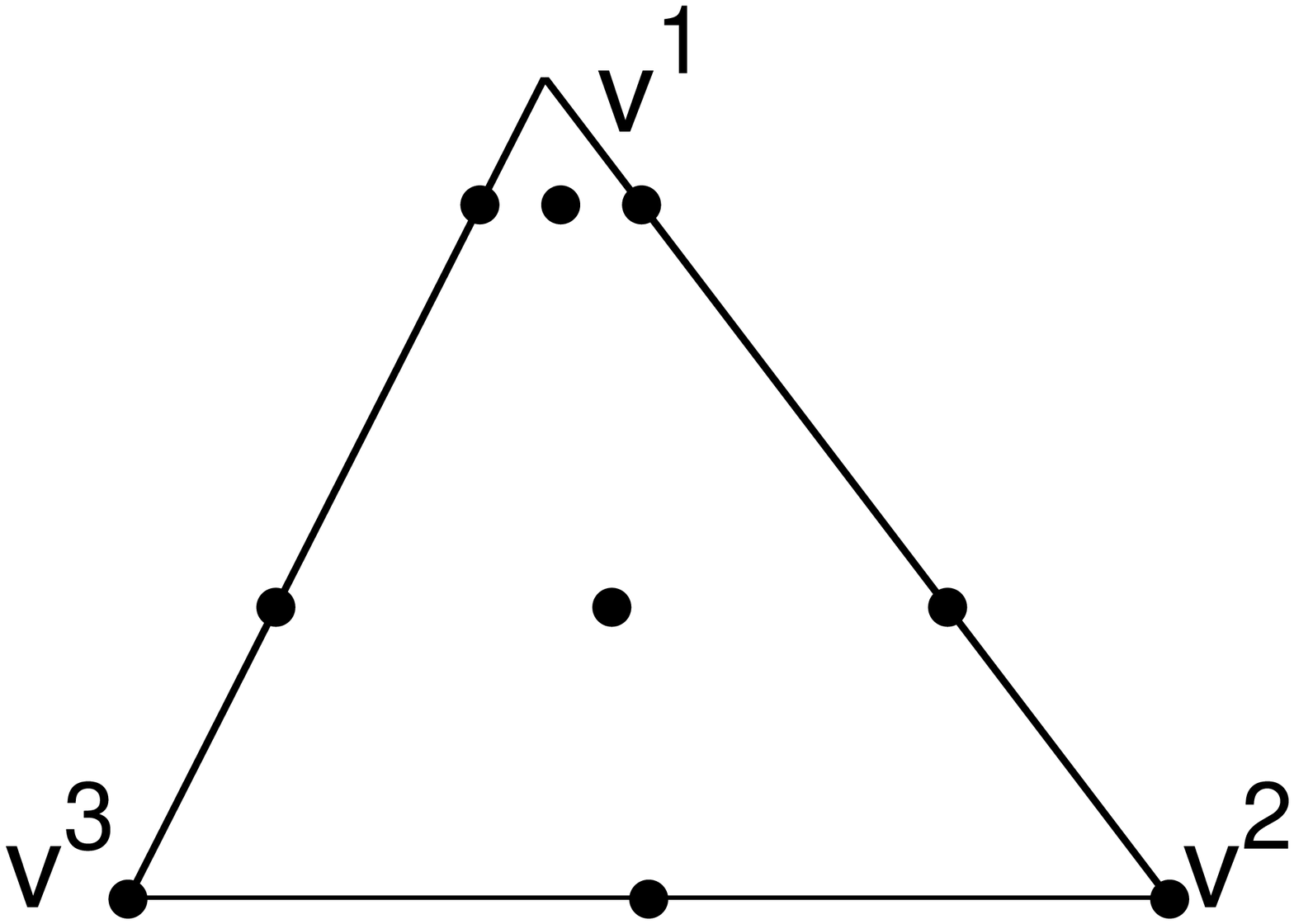}\label{fig:quadS2K_1}}
\subfigure[${\bf g}_2(S^2)$]{\includegraphics[width=0.24\textwidth]{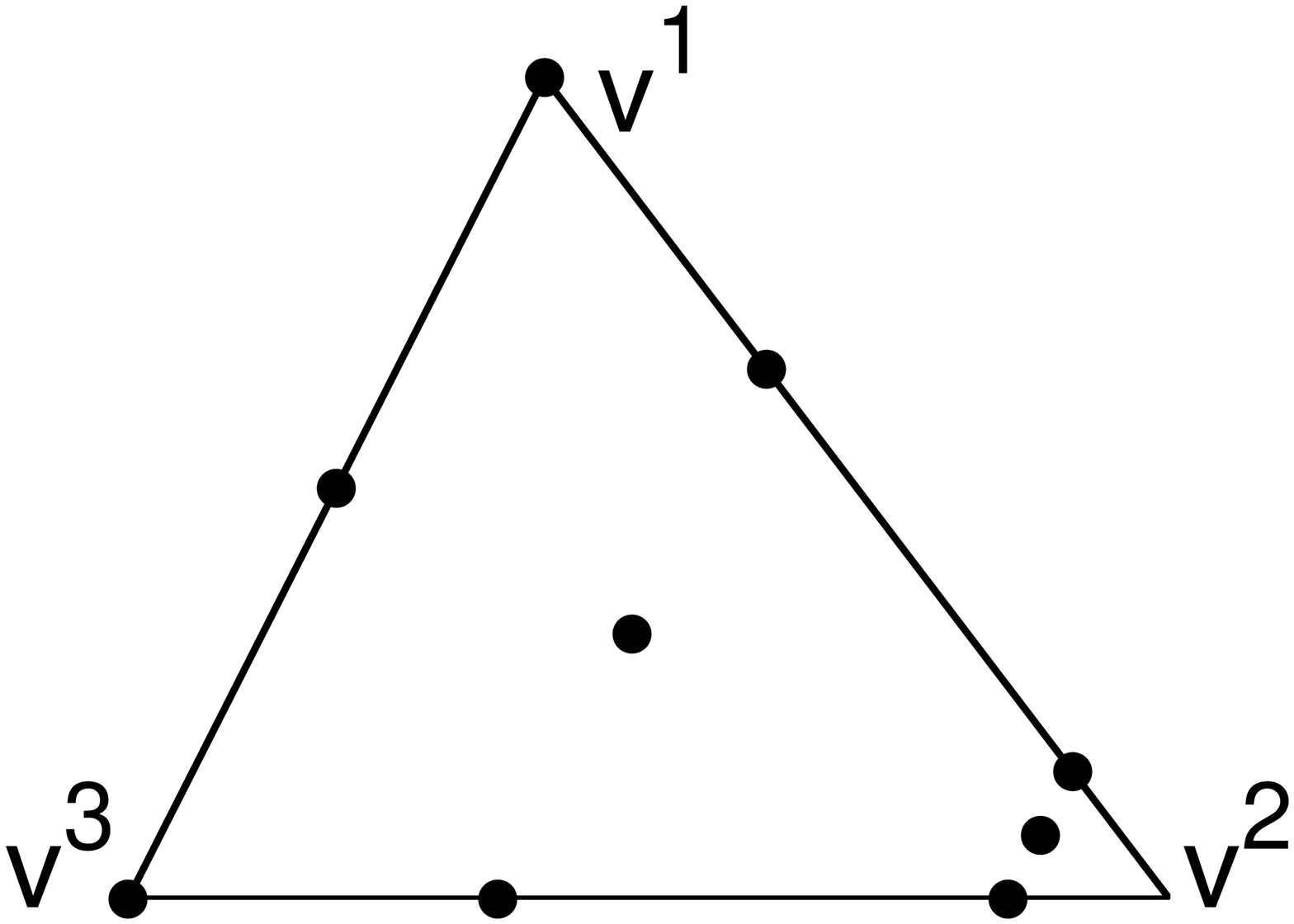}\label{fig:quadS2K_2}}
\subfigure[${\bf g}_3(S^2)$]{\includegraphics[width=0.24\textwidth]{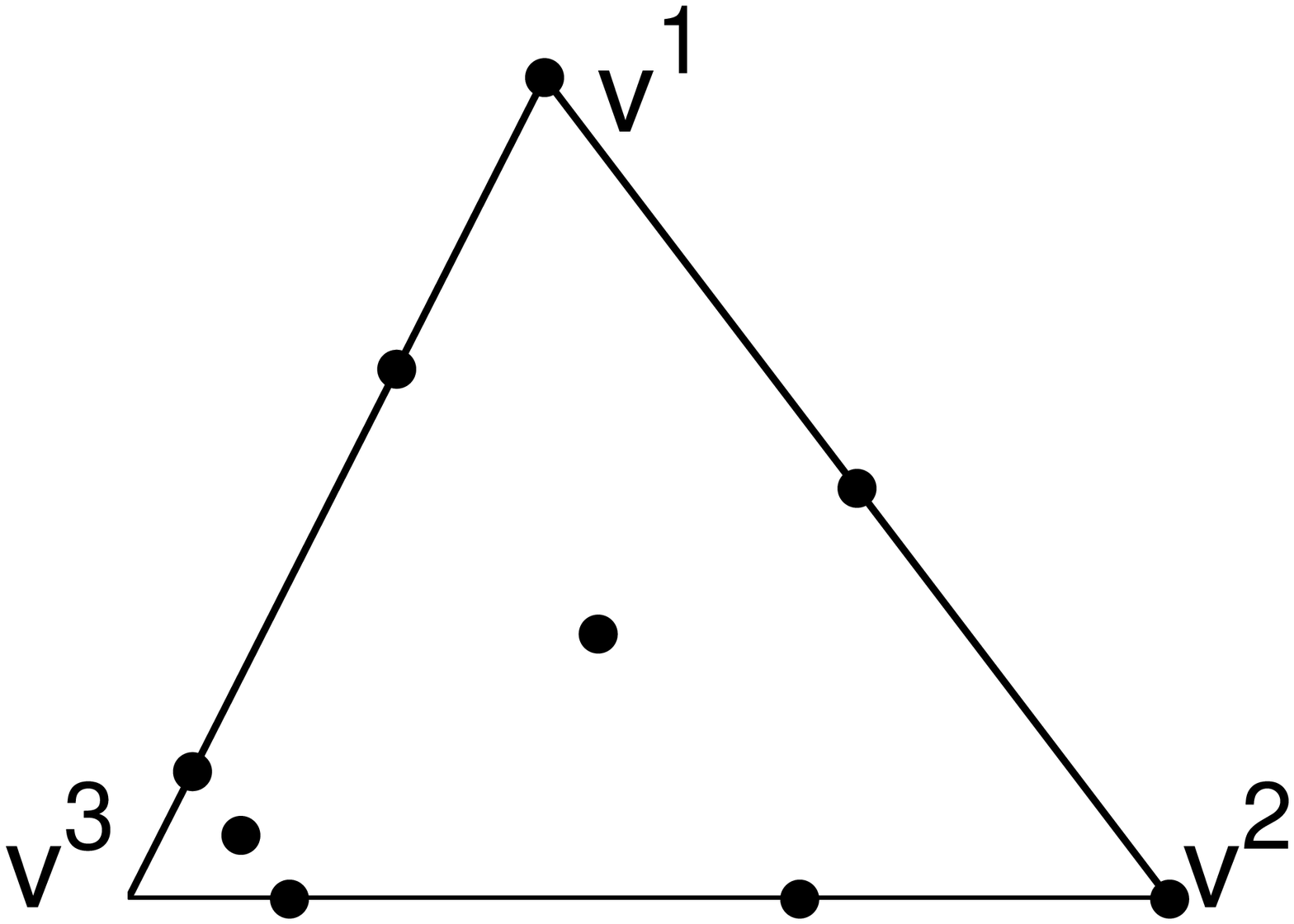}\label{fig:quadS2K_3}}
\caption{Illustration of quadrature points mapping from rectangle element to triangle element.}
\end{figure}

Let $\left\{\widehat{u}^\alpha: \alpha = 1, 2, 3\right\}$ denote the Gauss-Lobatto quadrature points on $\left[-\frac{1}{2},\frac{1}{2}\right]$ with weights $\widehat{w}_\alpha$ (in Table \ref{tab:quadG-L}), which is exact for one variable polynomial of degree 3. We have $\left\{v^\beta: \beta = 1, 2, 3\right\}$ (including the left boundary $v^1=-\frac{1}{2}$) denote the 3-point Gauss-Radau quadrature points on $\left[-\frac{1}{2},\frac{1}{2}\right]$ with weights $w_\beta$ (in Table \ref{tab:quadG-R}), which is exact for one variable polynomial of degree 4. For a two-variable polynomial $p(u,v)$, we use tensor product of 3-point Gauss-Lobatto for $u$ and 3-point Gauss Radau for $v$ as the quadrature rule on the square. The quadrature points can be written as $S^2=\left\{\left(\widehat{u}^\alpha,v^\beta\right): \alpha = 1, 2, 3; \beta = 1, 2,3\right\}$ with weights $\widehat{w}_\alpha w_\beta$, listed in Figure \ref{fig:quadS2K_0}.

\begin{table}[htbp]
   \centering
   \begin{minipage}{.35\textwidth}
   \centering
   \begin{tabular}{|c||c|c|c|}
      \hline
      $\alpha$ & 1 & 2 & 3 \\
      \hline
      \hline
      $\widehat{u}^\alpha$ & $-\frac{1}{2}$ & $0$ & $\frac{1}{2}$ \\
      \hline
      $\widehat{w}_\alpha$ & $\frac{1}{6}$ & $\frac{2}{3}$ & $\frac{1}{6}$\\
      \hline
   \end{tabular}
   \caption{3-point Gauss-Lobatto quadrature rule}
   \label{tab:quadG-L}
   \end{minipage}%
   \begin{minipage}{.6\textwidth}
   \centering
   \begin{tabular}{|c||c|c|c|}
      \hline
      $\beta$ & 1 & 2 & 3 \\
      \hline
      \hline
      $v^\beta$ & $-\frac{1}{2}$ & $\frac{1}{10}(1-\sqrt{6})$ & $\frac{1}{10}(1+\sqrt{6})$ \\
      \hline
      $w_\beta$ & $\frac{1}{9}$ & $\frac{1}{36}(16+\sqrt{6})$ & $\frac{1}{36}(16-\sqrt{6})$\\
      \hline
   \end{tabular}
   \caption{3-point Gauss-Radau quadrature rule}
   \label{tab:quadG-R}
\end{minipage}
\end{table}

Without loss of generality, we assume the orientation of the three vertices ${\bf v}^1$, ${\bf v}^2$ and ${\bf v}^3$ is marked clockwise. We define the following three functions as projections from the square to triangle $K$, mapping the top edge of the square into one vertex and the other three edges to the edges of $K$.
\begin{eqnarray}
{\bf g}_1(u,v) &=& \left( \frac{1}{2}+v\right){\bf v}^1 + \left( \frac{1}{2}+u\right)\left( \frac{1}{2}-v\right){\bf v}^2 + \left( \frac{1}{2}-u\right)\left( \frac{1}{2}-v\right){\bf v}^3, \nonumber \\
{\bf g}_2(u,v) &=& \left( \frac{1}{2}+v\right){\bf v}^2 + \left( \frac{1}{2}+u\right)\left( \frac{1}{2}-v\right){\bf v}^3 + \left( \frac{1}{2}-u\right)\left( \frac{1}{2}-v\right){\bf v}^1,\nonumber \\
{\bf g}_3(u,v) &=& \left( \frac{1}{2}+v\right){\bf v}^3 + \left( \frac{1}{2}+u\right)\left( \frac{1}{2}-v\right){\bf v}^1 + \left( \frac{1}{2}-u\right)\left( \frac{1}{2}-v\right){\bf v}^2. \nonumber
\end{eqnarray}
Under each projection ${\bf g}_i$ ($i = 1, 2$ or $3$), the quadrature points $S^2$ are mapped onto the triangle $K$, i.e. ${\bf g}_i(S^2)$, as in Figure \ref{fig:quadS2K_1} - Figure \ref{fig:quadS2K_3} . Let $S^2_K={\bf g}_1(S^2)\cup {\bf g}_2(S^2)\cup {\bf g}_3(S^2)$.

We use ${\bf g}_i (i=1,2,3)$ and $S^2$ to construct our triangle element quadrature rule. Let $p_K(x,y)$ be a two-variable polynomial of degree 2 with cell average $\overline{p_K}=\frac{1}{|K|}\int_K p_K(x,y) dxdy$, then we have,
\begin{eqnarray}\label{eq:uavecomb}
\overline{p_K} &=& \frac{1}{|K|} \int_K \! p_K(x,y) \, dA
= \frac{1}{|K|} \int_{-\frac{1}{2}}^{\frac{1}{2}} \int_{-\frac{1}{2}}^{\frac{1}{2}} \! p_K({\bf g}_i(u,v)) \left|\frac{\partial {\bf g}_i(u,v)}{\partial (u,v)}\right| \, dudv , \quad i=1,2,3 \nonumber \\
&=&  \int_{-\frac{1}{2}}^{\frac{1}{2}} \int_{-\frac{1}{2}}^{\frac{1}{2}} \! p_K({\bf g}_i(u,v))2(\frac{1}{2}-v) \, dudv
= \sum_{\alpha=1}^3\sum_{\beta=1}^3 p_K({\bf g}_i(\widehat{u}^\alpha,v^\beta))2(\frac{1}{2}-v^\beta) \widehat{w}_\alpha w_\beta  \nonumber \\
&=& \sum_{i=1}^3 \sum_{\alpha=1}^3\sum_{\beta=1}^3 p_K({\bf g}_i(\widehat{u}^\alpha,v^\beta))\frac{2}{3}(\frac{1}{2}-v^\beta) \widehat{w}_\alpha w_\beta
= \sum_{{\bf x} \in S^2_K} p_K({\bf x})w_{\bf x}.
\end{eqnarray}
With the three vertices ${\bf v}^1$, ${\bf v}^2$ and ${\bf v}^3$ orientated clockwise, we have the Jacobian $\left|\frac{\partial {\bf g}_i(u,v)}{\partial (u,v)}\right| = 2|K|(\frac{1}{2}-v)$. Notice that $p_K({\bf g}_i(u,v))2(\frac{1}{2}-v)$ is a polynomial of $u$ and $v$ with degree 2 and degree 3, therefore the quadrature rule on $S^2$ is exact.

It is easy to show the weights $w_{\bf x}$ for quadrature points ${\bf x} \in S^2_K$ are non-negative, then we can rewrite (\ref{eq:uavecomb}) as a combination of quadrature points, see below,
\begin{equation}
\overline{p_K} = \sum_{{\bf x} \in S^2_K\backslash \{{\bf v}^1, {\bf v}^2, {\bf v}^3\}} p_K({\bf x})w_{\bf x} + \sum_{i=1}^3 p_K({\bf v}^i)\overline{w}_i.
\end{equation}

We are interested in the weights $\left\{\overline{w}_i\right\}_{i=1}^3$ for all three vertices ${\bf v}^1$, ${\bf v}^2$ and ${\bf v}^3$. Let us take ${\bf v}^1$ for example. Notice that ${\bf g}_2(-\frac{1}{2},-\frac{1}{2})$ and ${\bf g}_3(\frac{1}{2},-\frac{1}{2})$ are the same point $(1, 0, 0)$, i.e., ${\bf v}^1$. Therefore, the weight of $(1, 0, 0)$ is
\begin{eqnarray}
\overline{w}_1 &=& \frac{2}{3}\left[\frac{1}{2}-\left(-\frac{1}{2}\right)\right]\widehat{w}_1 w_1 + \frac{2}{3}\left[\frac{1}{2}-\left(-\frac{1}{2}\right)\right]\widehat{w}_{3}w_{1}
= \frac{2}{3}\left(\widehat{w}_1+\widehat{w}_{3}\right)w_{1}
= \frac{2}{81} \label{eq:weight4vertex}.
\end{eqnarray}

\begin{rem}
This section only provide one way to construct a quadrature rule on any triangle, with three vertices and edge centers included as quadrature points. The goal is to show that one can find such a quadrature rule with positive weights for quadrature points.
\end{rem}

\bigskip
{\centering \bf{Acknowledgements.}}
Huang's work is supported by Natural Science Foundation of Zhejiang Province grant No.LY14A010002 and No.LY12A01009.

\bigskip
\newpage


\end{document}